\documentclass[final,3p,times]{elsarticle}
\usepackage{amssymb}
\usepackage{mathrsfs}
\usepackage{amsfonts}
\usepackage{amsmath,amssymb}
\usepackage[all]{xy}
\usepackage{latexsym}
\usepackage{amsthm,a4wide,color}
\usepackage{amsmath,amscd,verbatim,enumerate}
\usepackage{hyperref}
\usepackage{graphicx}
\usepackage{extarrows}
\usepackage{lineno}
\biboptions{sort&compress}
\usepackage{threeparttable}
\usepackage{psfrag, epsfig}
\usepackage{multirow}

\usepackage{mathrsfs}
\usepackage{longtable}
\usepackage{rotating}
\usepackage{float}
\usepackage{extarrows}
\usepackage{bbm}

\input diagxy

\theoremstyle{plain}
  \newtheorem{thm}{Theorem}[section]
   \newtheorem{theorem}[thm]{Theorem}
  \newtheorem{lemma}[thm]{Lemma}
  \newtheorem{proposition}[thm]{Proposition}
  
\theoremstyle{definition}
  \newtheorem{definition}[thm]{Definition}
  
  \newtheorem{example}[thm]{Example}
  \newtheorem{remark}[thm]{Remark}

\begin{document}

\newcommand{\oto}{{\to\hspace*{-3.1ex}{\circ}\hspace*{1.9ex}}}
\newcommand{\lam}{\lambda}
\newcommand{\da}{\downarrow}
\newcommand{\Da}{\Downarrow\!}
\newcommand{\D}{\Delta}
\newcommand{\ua}{\uparrow}
\newcommand{\ra}{\rightarrow}
\newcommand{\la}{\leftarrow}
\newcommand{\lra}{\longrightarrow}
\newcommand{\lla}{\longleftarrow}
\newcommand{\rat}{\!\rightarrowtail\!}
\newcommand{\up}{\upsilon}
\newcommand{\Up}{\Upsilon}
\newcommand{\ep}{\epsilon}
\newcommand{\ga}{\gamma}
\newcommand{\Ga}{\Gamma}
\newcommand{\Lam}{\Lambda}

\newcommand{\TF}{\mathbb{F}}
\newcommand{\TU}{\mathbb{U}}
\newcommand{\TW}{\mathbb{W}}
\newcommand{\TG}{\mathbb{G}}
\newcommand{\TB}{\mathbb{B}}
\newcommand{\TV}{\mathbb{V}}
\newcommand{\TT}{\mathcal{T}}
\newcommand{\TS}{\mathcal{S}}

\newcommand{\CU}{\mathscr{U}}
\newcommand{\CW}{\mathscr{W}}
\newcommand{\CCF}{\mathscr{F}}
\newcommand{\CCG}{\mathscr{G}}
\newcommand{\CVV}{\mathscr{V}}
\newcommand{\CN}{{\mathcal{N}}}

\newcommand{\SU}{{\mathcal{U}}}
\newcommand{\SV}{{\mathcal{V}}}
\newcommand{\SF}{{\mathcal{F}}}
\newcommand{\SG}{{\mathcal{G}}}
\newcommand{\SB}{{\mathcal{B}}}
\newcommand{\ST}{{\mathcal{T}}}

\newcommand{\DV}{{\mathbf{V}}}

\newcommand{\DB}{{\mathbf{B}}}

\numberwithin{equation}{section}
\renewcommand{\theequation}{\thesection.\arabic{equation}}

\begin{frontmatter}
\title{ A novel axiomatic approach to $L$-valued rough sets within an $L$-universe via inner product and outer product of $L$-subsets}


\author{Lingqiang Li\corref{cor}}
\ead{lilingqiang0614@126.com}\cortext[cor]{Corresponding author, Tel:
+86 15206506635}
\author{Qiu Jin}
\ead{jinqiu79@126.com}

\address{Department of Mathematics, Liaocheng University, Liaocheng 252059, P.R.China}

\begin{abstract}

The fuzzy rough approximation operator serves as the cornerstone of fuzzy rough set theory and its practical applications. Axiomatization is a crucial approach in the exploration of fuzzy rough sets, aiming to offer a clear and direct characterization of fuzzy rough approximation operators. Among the fundamental tools employed in this process, the inner product and outer product of fuzzy sets stand out as essential components in the axiomatization of fuzzy rough sets. In this paper, we will develop the axiomatization of a comprehensive fuzzy rough set theory, that is, the so-called $L$-valued rough sets  with an $L$-set serving as the foundational universe (referred to as the $L$-universe) for defining $L$-valued rough approximation operators, where $L$ typically denotes a GL-quantale. Firstly, we give the notions of  inner product and outer product of two $L$-subsets within an $L$-universe and examine their basic properties. It is shown that these notions are extensions of the corresponding notion of  fuzzy sets within  a classical universe. Secondly, leveraging the inner product and outer product of $L$-subsets, we respectively characterize $L$-valued upper and  lower rough approximation operators generated by general, reflexive, transitive, symmetric, Euclidean, and median $L$-value relations on $L$-universe as well as  their compositions. Finally, utilizing the provided axiomatic characterizations, we present the precise examples for the least and largest equivalent $L$-valued upper and  lower rough approximation operators. Notably,  many existing axiom characterizations of fuzzy rough sets within classical universe can be viewed as direct consequences of our findings.
\end{abstract}


\end{frontmatter}

\section{Introduction}

Rough set is a fundamental methodology for managing information in the presence of uncertainty and incompleteness \cite{ZP82}. Fuzzy rough set  is a meaningful extension of rough set, which can deal with more complex uncertain problems. Constructivism and axiomatization are two basic research approaches of (fuzzy) rough sets.  Constructive approaches involve initiating from a binary (fuzzy) relation, covering, or neighborhood, etc., to create a pair of rough approximation operators ({\it RAOs} for short). These operators are subsequently utilized to forge connections with different structures including topology, lattice, and matroid \cite{GD24, KH, GL23,FF22}, or to address practical issues like decision-making, attribute reduction, and classification clustering \cite{HN24,CA22,AT23,DZ24,HZH21}. Axiomatic approaches, on the other hand, involve characterizing {\it RAOs} through a set of axioms or a single axiom \cite{YY983,FF22,XW22}. Yao \cite{YY983} established the framework of axiomatic approach to rough sets, Liu  \cite{GL13}  provided an approach through a single axiom.  Morsi \cite{MN100} characterized the fuzzy {\it RAOs}  arise from fuzzy similarity relations, Wu \cite{WZ159,WZ16} extended these results to general fuzzy relations.

%


The original fuzzy rough set uses the unit interval $[0,1]$ as the membership range of fuzzy sets (called $[0,1]$-rough set). However, as Goguen pointed out, in some cases, it may not be appropriate to use linear order to express membership. Therefore, some lattice structures $L$ (especially the residuated lattice) are proposed to replace $[0,1]$ as the range of membership. The corresponding fuzzy rough sets and rough approximation operators are called $L$-rough set and $L$-{\it RAOs} in this paper. For $L$ a residuated lattice, Radzikowska-Kerre \cite{AM05} introduced a constructive approach to $L$-relation based $L$-{\it RAOs}, which was studied and generalized by many researchers. Ma-Hu \cite{ZM23} discussed the topological and lattice structures induced by $L$-{\it RAOs},  Sun-Shi \cite {SY23} gave the representations of $L$-{\it RAOs}, Pang-Mi-Yao \cite{BP19} explored  more $L$-{\it RAOs} via three new types of $L$-relations, and Zhao-Hu \cite{ZX15} developed the  variable precision $L$-{\it RAOs}. For the axiomatic approach, She-Wang \cite{SYH09} characterized the $L$-{\it RAOs} by a set of axioms, Bao-Yang-She \cite{Bao18} and Wang \cite{CY20} characterized them with a single axiom, and Jin-Li \cite{JIN232} further characterized the variable precision $L$-{\it RAOs}.
Both from constructive and axiomatic approaches, Li-Jin-Hu \cite{LQ17} discussed three kinds of $L$-covering-based $L$-{\it RAOs}, and Li-Yang-Qiao \cite{YB23} extended them to $\beta$-covering case.

In the aforementioned $L$-rough set theory, the basic universe is defined as a nonempty set $X$ (referred to as classical universe). Currently, there is a new trend in the study of fuzzy sets, that is, many traditional fuzzy concepts defined on universe $X$ are extended to $L$-set of $X$ (referred to as $L$-universe). Pu-Zhang \cite{Pu12} expanded the conventional $L$-relations on classical universe to $L$-valued relations on $L$-universe. Demirci \cite{MD22} introduced the notion of multi-valued topologies on $L$-universe. Fang-Yue \cite{JF20} discussed the function space on $L$-universe. Building upon the concept of $L$-valued relations, Li-Yue \cite{LF19} proposed a model of $L$-valued rough sets on $L$-universe, which serves as a more comprehensive framework for $L$-rough sets, extending beyond the traditional $L$-relations. Wei-Pang-Mi \cite{BP21} and Chen-Wei \cite{XW22} further developed the theory of $L$-valued rough sets and provided axiomatic characterizations of the theory. More recently, Temraz-Saady investigated the concept of variable precision $L$-valued rough sets on $L$-universe. Saady-Rashed-Temraz \cite {EI24} then presented the concept of $L$-valued covering-based rough sets in the context of the $L$-universe, along with their applications in MADM.

After the comprehensive review, we can succinctly summarize the evolution of rough set theory and its fuzzy extension, as depicted in Figure \ref{fig11} below.
 Herein, $L$-valued rough sets are rooted in the $L$-universe, whereas other theories adhere to the classical universe. Across these models, two foundational research approaches are consistently employed: construction and axiomatization.

\begin{figure}[!ht]
\centering
\includegraphics[width=14cm,height=5cm]{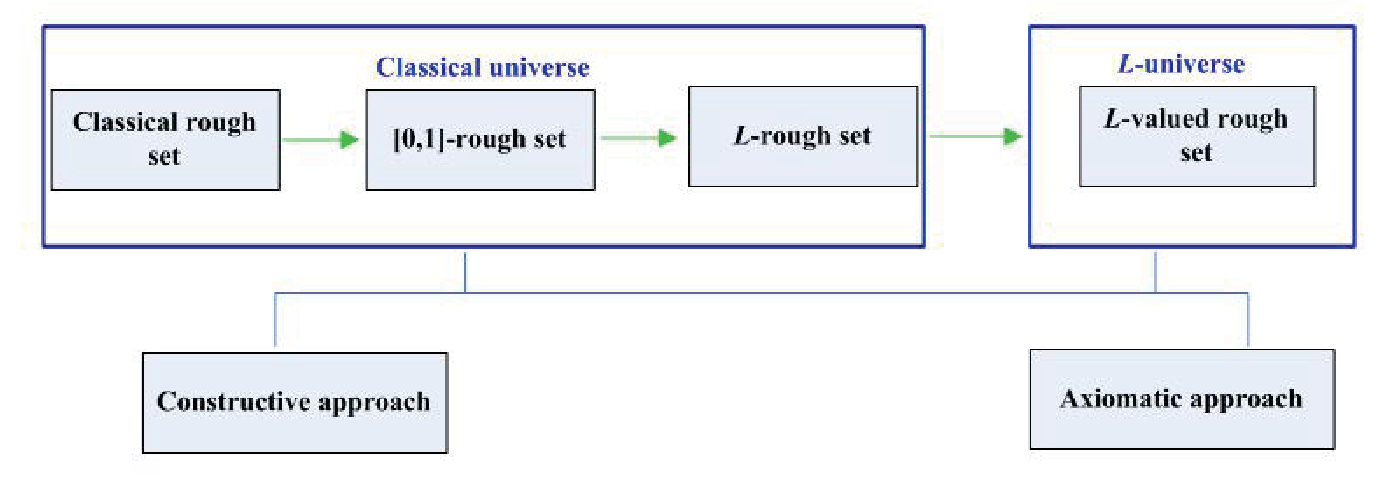}
\caption{ The evolution of rough set theory and its fuzzy extension} \label{fig11}
\end{figure}


Compared to the extensive research on $L$-rough sets based on classical universe, there has been relatively less exploration into $L$-valued rough sets on  $L$-universe. Thus, we shall delve deeper into this topic and provide a more comprehensive understanding. The primary motivations are as follows:

(1) The degree of intersection \cite{GG04} and the subsethood degree \cite{RB12} of two $L$-sets within a classical universe are crucial in $L$-rough sets and various other fuzzy structures. These concepts can facilitate the provision of natural explanations and concise characterizations for certain fuzzy notions \cite{GD24,EI,YB23,JS23}. To date, these two notions have not been presented within the context of an $L$-universe. Consequently, to enhance the fuzzy theory grounded in an $L$-universe, it is essential to introduce the degree of intersection and the subsethood degree for two $L$-subsets within the $L$-universe framework.

(2) Drawing from the research on $L$-rough sets within the classical universe as detailed in the references \cite{Bao18,JIN232,JJ23}, scholars have utilized the inner product (that is, the degree of intersection) and the outer product (defined through the subsethood degree) of the two $L$-sets to derive numerous interesting axiomatic characterizations on  $L$-{\it RAOs} that arise from general, reflexive, symmetric, transitive, Euclidean, and mediate $L$-relations, as well as their combinations. Nevertheless, the corresponding axiomatic characterizations on $L$-valued rough sets within the $L$-universe, as presented in \cite{LF19}, have not been established. It is thus inherently significant to address this gap.

To facilitate reading, we provide some symbols and abbreviations for commonly used terms in Table \ref{JX2}.

\begin{table}[!h] \footnotesize
\centering \label{JX2}\caption{The symbols and abbreviations for commonly used terms}
\setlength{\tabcolsep}{2mm}{
\begin{tabular}{llll}
\hline
Symbol& Meaning\\
\hline
$X$&  \ \ \ \ \ \ \ A nonempty set, called a classical universe\\
$\mathbb{U}$&  \ \ \ \ \ \ \ A fixed $L$-set in $X$, called an $L$-universe\\
$P(\mathbb{U})$&  \ \ \ \ \ \ \ All $L$-subsets in $\mathbb{U}$\\
\hline
Full name& Abbreviation\\
\hline
$L$-valued relation on $L$-universe&  \ \ \ \ \ \ \ {\it Lvr}\\
$L$-valued lower  rough approximation operator  based on {\it Lvr}&  \ \ \ \ \ \ \ {\it LVLRAO}\\
$L$-valued upper  rough approximation operator  based on {\it Lvr}&  \ \ \ \ \ \ \ {\it LVHRAO}\\
\hline
\end{tabular}}\label{JX2}
\end{table}

The paper is organized as follows: Section 2 lays down the fundamental concepts and notation. Section 3 introduces the concept of inner product (that is, the degree of intersection) of the two $L$-subsets in $L$-universe and employs this to provide axiomatic characterizations on 18 types of {\it LVHRAO} generated by general, reflexive, transitive, symmetric, Euclidean, median {\it Lvr}s and their compositions. Section 4 defines the notion  of the outer product (and subsethood degree)  of the two $L$-subsets in $L$-universe and utilizes this to offer axiomatic characterizations on 18 types of {\it LVLRAO} generated by general, reflexive, transitive, symmetric, Euclidean, median {\it Lvr}s and their compositions. Employing the specified axiomatic descriptions, we provide the precise examples for the least and largest equivalent {\it LVHRAO} and {\it LVLRAO}. Section 5 concludes and closes the paper.

The mind map of this paper is summarized in Figure \ref{fig44} below.

\begin{figure}[!ht]
\centering
\includegraphics[width=14cm,height=4cm]{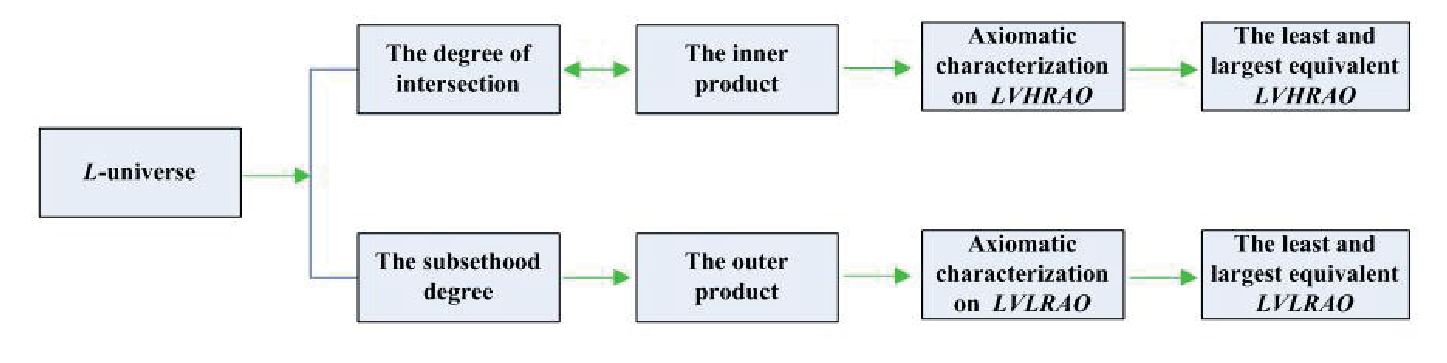}
\caption{The mind map of this paper} \label{fig44}
\end{figure}

%
%
%

\section{Preliminaries}

This section shall review some knowledge about fuzzy rough sets as preliminary.

\subsection{The $L$-rough sets on classical universal and their three axiomatic approaches}

Let $(L,\circledast)$  be a  {\it complete residuated lattice}, that is, $L=(L,\wedge,\vee, 1,0)$  is a complete lattice and $\circledast:L\times L\longrightarrow L$  a commutative, associative binary operation fulfilling

(i) $\forall \alpha\in L, \alpha\circledast 1=1$,

(ii) $\forall \alpha, \beta_j(j\in J) \in L$,  $\alpha\circledast \bigvee\limits_{j}\beta_j= \bigvee\limits_{j}(\alpha\circledast \beta_j)$.

A binary operation  $\rightarrow:L\times L\longrightarrow L$ determined by $\forall \alpha, \beta\in L, \alpha \rightarrow \beta=\bigvee\{\gamma\in L:\alpha \circledast \gamma \leq \beta\}$ is termed the {\it residuated implication} of $\circledast$.

Some notable  properties of $\circledast, \rightarrow$  are gathered below.

\begin{lemma} \cite{R90} \label{jb00} Suppose $\alpha, \beta, \gamma, \alpha_{i}(i\in \Lambda) \in L$.

{\rm (1)} $\alpha\circledast \gamma\leq \beta \Leftrightarrow \gamma\leq \alpha\rightarrow \beta$.

{\rm (2)} $(\alpha\rightarrow \beta)\circledast \alpha\leq \beta$.

{\rm (3)} $\alpha\rightarrow (\beta\rightarrow \gamma)=(\alpha\circledast \beta)\rightarrow \gamma$.

{\rm (4)} $\alpha\leq \beta \Leftrightarrow \alpha\rightarrow \beta=1$.

{\rm (5)} $\beta\rightarrow (\bigwedge\limits_i \alpha_{i})=\bigwedge\limits_i(\beta\rightarrow \alpha_{i})$.

{\rm (6)} $(\bigvee\limits_i \alpha_{i})\rightarrow
\beta=\bigwedge\limits_i(\alpha_{i}\rightarrow \beta)$.
\end{lemma}

Let $X$ denote a nonempty set and $2^X$ represent the power set of $X$. A mapping   $W: X\longrightarrow L$ is termed an $L$-set on $X$.  All $L$-sets on $X$ are recorded  as $L^X$. Let $W,V \in L^X$  and $\diamond\in \{\vee,\wedge, \circledast, \rightarrow \}$. One define $\sim W, W\diamond V: X\longrightarrow L$ through $\forall d\in X$,
$$\sim W(d)=W(d)\rightarrow 0, (W\diamond V)(d)=W(d)\diamond V(d).$$
In addition, for $\alpha\in L$, we also use $\alpha$ to denote the constant $L$-set on $X$ valued $\alpha$.

For any $M\in 2^X$, we  define $\chi_M\in L^X$ by $\chi_M(x)=\left\{\begin{array}{ll}1,& x\in M,\\ 0,&x\not\in M.\end{array}\right.$

\begin{definition} \cite{JIN232,AM05,BP19} One mapping $\mathfrak{R}: X\times X\longrightarrow L$ is named an $L$-relation on $X$. The couple $(X, \mathfrak{R})$ is named an $L$-approximation space. Furthermore,  $\mathfrak{R}$ is referred:

(1) reflexive provided  $\mathfrak{R}(a,a)=1$;

(2) symmetric whenever  $\mathfrak{R}(a,d)=\mathfrak{R}(d,a)$;

(3) transitive provided  $\mathfrak{R}(a,d)\circledast \mathfrak{R}(d,h)\leq \mathfrak{R}(a,h)$;

\noindent for any $a,d,h\in X$. $\mathfrak{R}$ is called equivalent if (1)-(3) are satisfied.  \end{definition}

\begin{definition} \cite{AM05} \label{deffr} Suppose $(X, \mathfrak{R})$ is an $L$-approximation space and $\underline{\mathfrak{R}}, \overline{\mathfrak{R}}: L^X\longrightarrow L^X$ are mappings determined by $\forall Q\in
L^X, o\in X$:
$$\underline{\mathfrak{R}}  (Q)(o)=\bigwedge\limits_{d\in X}\Big(\mathfrak{R}(d,o)\rightarrow Q(d)\Big), \overline{\mathfrak{R}} (Q)(o)=\bigvee\limits_{d\in X}\Big(\mathfrak{R}(d,o)\circledast Q(d)\Big).$$ Then the couple  $(\underline{\mathfrak{R}}  (Q), \overline{\mathfrak{R}}(Q))$ is termed  $L$-rough set about $Q$. Furthermore, $\underline{\mathfrak{R}}$  is termed
$L$-lower  rough approximation operator, and  $\overline{\mathfrak{R}} (Q)$ is termed
$L$-upper  rough approximation operator.
\end{definition}

In the literatures, three axiomatic approaches to $L$-rough sets were presented.

\begin{proposition} \label{fuzzy1} (Approach 1, \cite{SYH09})  Let $\mathbf{H}, \mathbf{L}: L^X\longrightarrow L^X$ be mappings.

(1) There exists an unique $L$-relation $\mathfrak{R}$ on $X$ s.t. $\mathbf{H}=\overline{\mathfrak{R}}$ iff (A1) and (A2) hold. Furthermore,  $\mathfrak{R}$ is equivalent  iff (A1)-(A4) hold. Where, $\forall W, W_i(i\in I)\in L^X, \forall \alpha\in L$:

{(A1)}  $\mathbf{H}   (\alpha\circledast W)=\alpha\circledast \mathbf{H}(W)$.

{(A2)}  $\mathbf{H}   (\bigvee\limits_{i}W_i)=\bigvee\limits_{i} \mathbf{H}   (W_i)$.

{(A3)} $\mathbf{H}(W)\geq W$.

{(A4)} $\mathbf{H}\mathbf{H}(W)\leq \mathbf{H}(W)$.

(A5) $\forall c,d\in X, \mathbf{H}(\chi_{\{d\}})(c)=\mathbf{H}(\chi_{\{c\}})(d)$.

(2) There exists an unique $L$-relation $\mathfrak{R}$ on $X$ s.t. $\mathbf{L}=\underline{\mathfrak{R}}$ iff (B1) and (B2) hold. Furthermore,  $\mathfrak{R}$ is equivalent  iff (B1)-(B5) hold. Where, $\forall W, W_i(i\in I)\in L^X, \forall \alpha\in L$:

(B1)  $\mathbf{L}(\alpha\rightarrow W)=\alpha\rightarrow\mathbf{L}(W)$.

(B2)  $\mathbf{L}(\bigwedge\limits_{i}W_i)=\bigwedge\limits_{i} \mathbf{L}   (W_i)$.

(B3) $\mathbf{L}(W)\leq W$.

(B4)   $\mathbf{L}\mathbf{L}(W)\geq \mathbf{L}(W)$.

(B5)   $\forall c, d\in X, \mathbf{H}(\chi_{\{d\}}\rightarrow \alpha )(c)=\mathbf{H}(\chi_{\{c\}}\rightarrow \alpha)(d)$.

\end{proposition}

\begin{proposition} \label{fuzzy2} (Approach 2, \cite{Bao18})  Let $\mathbf{H}, \mathbf{L}: L^X\longrightarrow L^X$ be mappings.

(1) There exists an unique $L$-relation $\mathfrak{R}$ on $X$ s.t. $\mathbf{H}=\overline{\mathfrak{R}}$ iff $\mathbf{H}$ fulfills the single axiom: $\forall W_i(i\in I)\in L^X, \forall \alpha_i(i\in I)\in L$, $\mathbf{H}(\bigvee \limits_{i}(\alpha_i\circledast W_i))=\bigvee\limits_{i} (\alpha_i\circledast\mathbf{H}(W_i)).$

(2) There exists an uniquely equivalent $L$-relation $\mathfrak{R}$ on $X$ s.t. $\mathbf{H}=\overline{\mathfrak{R}}$ iff $\mathbf{H}$ fulfills the single axiom: $\forall W, W_i(i\in I)\in L^X, \forall \alpha, \alpha_i(i\in I)\in L$,

$ (\alpha\circledast \mathbf{H}^{-1}(W)) \vee \mathbf{H}(\bigvee\limits_{i} (\alpha_i\circledast W_i))=(\alpha\circledast \mathbf{H}(W))\vee (\bigvee\limits_i (\alpha_i\circledast W_i))\vee (\bigvee\limits_{i} (\alpha_i\circledast \mathbf{H}(W_i)))\vee (\bigvee\limits_{i}(\alpha_i\circledast\mathbf{H}\mathbf{H}(W_i))).$

(3) There exists an unique $L$-relation $\mathfrak{R}$ on $X$ s.t. $\mathbf{L}=\underline{\mathfrak{R}}$ iff $\mathbf{L}$ fulfills the single axiom: $\forall W_i(i\in I)\in L^X, \forall \alpha_i(i\in I)\in L$. $\mathbf{L}(\bigwedge \limits_{i}(\alpha_i\rightarrow W_i))=\bigwedge\limits_{i} (\alpha_i\rightarrow\mathbf{L}(W_i)).$

(4) There exists an uniquely equivalent $L$-relation $\mathfrak{R}$ on $X$ s.t. $\mathbf{L}=\underline{\mathfrak{R}}$ iff $\mathbf{L}$ fulfills the single axiom: $\forall W, W_i(i\in I)\in L^X, \forall \alpha, \alpha_i(i\in I)\in L$,

$ (\alpha\rightarrow \mathbf{L}_{-1}(W)) \wedge \mathbf{L}(\bigwedge\limits_{i} (\alpha_i\rightarrow W_i))=(\alpha\rightarrow \mathbf{L}(W))\wedge (\bigwedge\limits_i (\alpha_i\rightarrow W_i))\wedge (\bigwedge\limits_{i} (\alpha_i\rightarrow \mathbf{L}(W_i)))\wedge (\bigwedge\limits_{i}(\alpha_i\rightarrow\mathbf{L}\mathbf{L}(W_i))).$

\end{proposition}

\begin{definition} \label{defset} \cite{RB12,Bao18} (1)
Define mapping  $\mathbf{I}: L^X\times L^X\longrightarrow L$  by $\forall M, Q\in L^X$,
$\mathbf{I}(M,Q)=\bigvee\limits_{d\in X}\Big(M(d)\circledast Q(d)\Big).$  Then  $ \mathbf{I}(M,Q)$ is named the  inner product (i.e., the degree of intersection) of $M,Q$.

 (2) Define mapping  $\mathbf{S}: L^X\times L^X\longrightarrow L$  by $\forall M, Q\in L^X$,
$\mathbf{S}(M,Q)=\bigwedge\limits_{d\in X}\Big(M(d)\rightarrow Q(d)\Big).$  Then  $ \mathbf{S}(M,Q)$ is named the  subset degree of $M,Q$.

 (3) Define mapping  $\mathbf{O}: L^X\times L^X\longrightarrow L$  by $\forall M, Q\in L^X$,
$\mathbf{O}(M,Q)=\mathbf{S}(\sim M,Q)=\bigwedge\limits_{d\in X}\Big(\sim M(d)\rightarrow Q(d)\Big).$  Then  $\mathbf{O}(M,Q)$ is named the  outer product of $M,Q$.
\end{definition}

\begin{definition}\cite{Bao18} Let $\mathbf{H}, \mathbf{L}: L^X\longrightarrow L^X$ be mappings.

(1) The mapping  $\mathbf{H}^{-1}: L^X\longrightarrow   L^X$ determined by $\forall Q\in  L^X, \forall d\in X$, $ \mathbf{H}^{-1}(Q)(d)=\mathbf{I}\Big(\mathbf{H}(\chi_{\{d\}}), Q\Big)$ is called the upper inverse mapping of $\mathbf{H}$.

(2) The mapping  $\mathbf{L}_{-1}: L^X\longrightarrow   L^X$ determined by $\forall Q\in  L^X, \forall d\in X$, $ \mathbf{L}_{-1}(Q)(d)=\mathbf{O}\Big(\mathbf{L}(\chi_{X-\{d\}}), Q\Big)$ is called the lower inverse mapping of $\mathbf{L}$.
\end{definition}

\begin{proposition} \label{fuzzy3} (Approach 3, \cite{Bao18})  Let $\mathbf{H}, \mathbf{L}: L^X\longrightarrow L^X$ be mappings.

(1) There exists an unique $L$-relation $\mathfrak{R}$ on $X$ s.t. $\mathbf{H}=\overline{\mathfrak{R}}$ iff $\mathbf{H}$ fulfills the single axiom
$\forall M,Q\in  L^X$, $\mathbf{I}(M, \mathbf{H}(Q))=\mathbf{I}(Q, \mathbf{H}^{-1}(M)).$

(2) There exists an uniquely equivalent $L$-relation $\mathfrak{R}$ on $X$ s.t. $\mathbf{H}=\overline{\mathfrak{R}}$ iff $\mathbf{H}$ fulfills the single axiom
$\forall M,Q\in  L^X$, $\mathbf{I}(M, Q\vee \mathbf{H}(Q)\vee \mathbf{H}\mathbf{H}(Q))=\mathbf{I}(Q, \mathbf{H}(M)).$

(3) There exists an unique $L$-relation $\mathfrak{R}$ on $X$ s.t. $\mathbf{L}=\underline{\mathfrak{R}}$ iff $\mathbf{L}$ fulfills the single axiom
$\forall M,Q\in  L^X$, $\mathbf{O}(M, \mathbf{L}(Q))=\mathbf{O}(Q, \mathbf{L}_{-1}(M)).$

(4) There exists an uniquely equivalent $L$-relation $\mathfrak{R}$ on $X$ s.t. $\mathbf{L}=\underline{\mathfrak{R}}$ iff $\mathbf{L}$ fulfills the single axiom
$\forall M,Q\in  L^X$, $\mathbf{O}(M, Q\wedge\mathbf{L}(Q)\wedge \mathbf{L}\mathbf{L}(Q))=\mathbf{O}(Q, \mathbf{L}(M)).$
\end{proposition}

\begin{remark} Noticed that the characterizations provided by Approaches 2 and 3 only contain one axiom, so they are called single axiom characterizations in existing literature \cite{Bao18,JIN232}. Furthermore, Approach 3 exhibits greater conciseness compared to Approach 2. \end{remark}

\subsection{The $L$-valued rough sets on $L$-universal and two  axiomatic approaches}

A complete residuated lattice $(L,\circledast)$ is named a {\it GL-quantale} if it fulfills moreover $\forall \alpha,\beta\in L$, $\alpha\wedge \beta=\alpha\circledast(\alpha\rightarrow \beta)$.

Unless further specified, we always assume that $(L,\circledast)$  is a GL-quantale, its  further  properties  are gathered below.

\begin{lemma} \cite{R90} \label{jb} Suppose  $(L,\circledast)$  is a GL-quantale. Then $\forall \alpha, \beta, \gamma\in L$

(1)   $\beta\leq \gamma$ implies $\beta=\gamma\circledast (\gamma\rightarrow \beta)$.

(2)  $\alpha, \beta\leq \gamma$ implies $\alpha\circledast (\gamma\rightarrow \beta)=\beta\circledast (\gamma\rightarrow \alpha)$.

(3)  $\beta\leq \gamma$ implies $\gamma\circledast (\beta\rightarrow \alpha)=\gamma\wedge ((\gamma\rightarrow \beta)\rightarrow \alpha)$.
\end{lemma}

\begin{definition} We always let $\mathbb{U}: X\longrightarrow L$ be a fixed $L$-set of $X$, and call it an {\it $L$-universe} on $X$. Then a mapping $W: X\longrightarrow L$ is referred an $L$-subset in $\mathbb{U}$ if $W(z)\leq \mathbb{U}(z)$ for any $z\in X$.  All $L$-subsets in $\mathbb{U}$ are recorded  as $P(\mathbb{U})$ named the $L$-powerset of $\mathbb{U}$. Obviously, when $\mathbb{U}=1$, then $P(\mathbb{U})=L^X$.\end{definition}

 For each $d\in X$, the mappings $\mathbb{U}_{\{d\}}, \mathbb{U}_{X-{\{d\}}}:X\longrightarrow L$ belong to $P(\mathbb{U})$, where
$$\mathbb{U}_{\{d\}}(x)=\left\{\begin{array}{ll}\mathbb{U}(x),& x=d,\\ 0,&otherwise.\end{array}\right.,\ \mathbb{U}_{X-{\{d\}}}(x)=\left\{\begin{array}{ll}0,& x=d,\\ \mathbb{U}(x),&otherwise.\end{array}\right.$$

Note that $P(\mathbb{U})\subseteq L^X$, so for $W,V \in P(\mathbb{U})$  and $\diamond\in \{\vee,\wedge, \circledast, \rightarrow \}$,  one can define $W\diamond V$ as that in $L^X$. When $\diamond\in \{\vee,\wedge, \circledast \}$ then $W\diamond V \in P(\mathbb{U})$, but $W\rightarrow V, \sim W\not\in P(\mathbb{U})$ generally. Hence, we define $\neg W\in L^X$ by $\forall d\in X$, $\neg W(d)=\mathbb{U}(d)\circledast (W(d)\rightarrow 0)$, then $\neg W\in P(\mathbb{U})$.
%

\begin{definition} \cite{XW22,LF19} One mapping $\mathbb{R}: X\times X\longrightarrow L$ is named an $L$-valued relation ({\it Lvr} for short) on $\mathbb{U}$ provided $\forall a, d\in X, \mathbb{R}(a,d)\leq \mathbb{U}(a)\wedge \mathbb{U}(d)$. The couple $(\mathbb{U}, \mathbb{R})$ is named an $L$-valued approximation space ({\it LVAP}, for short). Furthermore,  $\mathbb{R}$ is referred:

(1) reflexive provided  $\mathbb{U}(a)\leq \mathbb{R}(a,a)$;

(2) symmetric whenever  $\mathbb{R}(a,d)=\mathbb{R}(d,a)$;

(3) transitive provided  $\mathbb{R}(a,d)\circledast (\mathbb{U}(d)\rightarrow \mathbb{R}(d,h))\leq \mathbb{R}(a,h)$;

\noindent for any $a,d,h\in X$.

$\mathbb{R}$ is named tolerance whenever it meets (1) + (2),  preordered whenever it meets (1) + (3), and equivalent  whenever it meets (1)-(3).
\end{definition}

\begin{definition} \cite{LF19} \label{defgraun} Suppose $(\mathbb{U}, \mathbb{R})$ is an {\it LVAP} and $\underline{\mathbb{R}}, \overline{\mathbb{R}}: P(\mathbb{U})\longrightarrow P(\mathbb{U})$ are mappings determined by $\forall Q\in
P(\mathbb{U}), o\in X$:
$$\underline{\mathbb{R}}  (Q)(o)=\bigwedge\limits_{d\in X}\Big(\mathbb{U}(o)\circledast \big(\mathbb{R}(d,o)\rightarrow Q(d)\big)\Big), \overline{\mathbb{R}} (Q)(o)=\bigvee\limits_{d\in X}\Big(\mathbb{R}(d,o)\circledast\big(\mathbb{U}(d)\rightarrow Q(d)\big)\Big).$$ Then the couple  $(\underline{\mathbb{R}}  (Q), \overline{\mathbb{R}}(Q))$ is termed  $L$-valued rough set about $Q$. Furthermore, $\underline{\mathbb{R}}$  is termed
$L$-valued lower  rough approximation operator ({\it LVLRAO} for short), and  $\overline{\mathbb{R}} (Q)$ is termed
$L$-valued  upper  rough approximation operator ({\it LVHRAO} for short).
\end{definition}

\begin{remark}  When $\mathbb{U}=1$, an {\it Lvr} $\mathbb{R}$ on $\mathbb{U}$ is precisely an $L$-relation on $X$, and $\forall Q\in
P(\mathbb{U}), a\in X$: $$\underline{\mathbb{R}}  (Q)(a)=\bigwedge\limits_{d\in X}\Big(\mathbb{R}(d,a)\rightarrow Q(d)\Big), \overline{\mathbb{R}} (Q)(a)=\bigvee\limits_{d\in X}\Big(Q(d)\circledast\mathbb{R}(d,a)\Big).$$ Hence   $(\underline{\mathbb{R}}  (Q), \overline{\mathbb{R}}(Q))$ becomes the $L$-rough set in Definition \ref{deffr}. So, $L$-valued rough set is a generalization of $L$-rough set \cite{AM05}.
\end{remark}

%
%
%
%
%



\begin{proposition} \label{theoremL1} (Approach 1, \cite{LF19})  Let $\mathbb{H}, \mathbb{L}: P(\mathbb{U})\longrightarrow P(\mathbb{U})$ be mappings.

{\rm (1)} There exists an unique {\it Lvr} $\mathbb{R}$ on $\mathbb{U}$ s.t. $\mathbb{H}=\overline{\mathbb{R}}$ iff (C1) and (C2) hold. Furthermore,  $\mathbb{R}$ is equivalent  iff (C1)-(C5) hold.

(C1)  $\mathbb{H}   (\beta \circledast (\alpha\rightarrow Q))=\beta \circledast(\alpha\rightarrow \mathbb{H} (Q))$ for any $\beta ,\alpha\in L, Q \in P(\mathbb{U})$ with  $\beta \leq \alpha$ and $\bigvee\limits_{d\in X}Q(d)\leq \alpha$.

(C2)   $\mathbb{H}   (\bigvee\limits_{i}W_i)=\bigvee\limits_{i} \mathbb{H}   (W_i)$ for any $W_i(i\in I)\subseteq P(\mathbb{U})$.

(C3)  $W\leq \mathbb{H}(W)$ for any $W \in P(\mathbb{U})$.

(C4)  $\mathbb{H}\mathbb{H}(W)\leq \mathbb{H}(W)$ for any $W \in P(\mathbb{U})$.

(C5) $\mathbb{H}(\mathbb{U}_{\{d\}})(h)=\mathbb{H}(\mathbb{U}_{\{h\}})(d)$ for any $d,h\in X$.

{\rm (2)} There exists an unique {\it Lvr} $\mathbb{R}$ on $\mathbb{U}$ s.t. $\mathbb{L}=\underline{\mathbb{R}}$ iff (D1) and (D2) hold. Furthermore,  $\mathbb{R}$ is equivalent  iff (D1)-(D5) hold.

(D1)  $\mathbb{L}(\mathbb{U}\wedge (\alpha\rightarrow W))=\mathbb{U}\wedge (\alpha\rightarrow \mathbb{L} (W))$ for any $\alpha\in L, W \in P(\mathbb{U})$.

(D2)  $\mathbb{L}   (\bigwedge\limits_{i}W_i)=\bigwedge\limits_{i} \mathbb{L}   (W_i)$ for any $W_i(i\in I)\subseteq P(\mathbb{U})$.

(D3)  $Q\geq \mathbb{L}(Q), \forall Q \in P(\mathbb{U})$.

(D4)  $\mathbb{L}\mathbb{L}(Q)\geq \mathbb{L}(Q), \forall Q \in P(\mathbb{U})$.

(D5) $\neg \mathbb{L}(\mathbb{U}_{X-\{d\}})(h)=\neg \mathbb{L}(\mathbb{U}_{X-\{h\}})(d)$ for any $d,h\in X$.
\end{proposition}

%
%
\begin{proposition} \label{WSC} (Approach 2, \cite{BP21})  Let $\mathbb{H},\mathbb{L}: P(\mathbb{U})\longrightarrow P(\mathbb{U})$ be mappings.

(1) There exists an unique {\it Lvr} $\mathbb{R}$ on $\mathbb{U}$ s.t. $\mathbb{H}=\overline{\mathbb{R}}$ iff $\mathbb{H}$ fulfills the single axiom: $\forall \Phi: P(\mathbb{U})\longrightarrow L$, $\forall \alpha,\beta\in L$ with $\alpha\leq \beta$, $\bigvee\limits_{d\in X}\sup \Phi(d)\leq \beta$,

$\mathbb{H}(\alpha\circledast(\beta\rightarrow \sup\Phi))=\alpha\circledast (\beta\rightarrow \bigvee\limits_{V\in P(\mathbb{U})}\Phi(V)\circledast \mathbb{H}(V)).$

(2) There exists an uniquely equivalent {\it Lvr} $\mathbb{R}$ on $\mathbb{U}$ s.t. $\mathbb{H}=\overline{\mathbb{R}}$ iff $\mathbb{H}$ fulfills the single axiom: $\forall \Phi: P(\mathbb{U})\longrightarrow L$, $\forall \alpha,\beta, \gamma\in L$ with $\alpha\leq \beta$ and $\bigvee\limits_{d\in X}\sup \Phi(d)\leq \beta$,

 $\mathbb{H}\Big(\gamma\circledast \bigvee\limits_{Q\in P(\mathbb{U})}\Phi(Q)\circledast(\neg \mathbb{H}(\neg Q)\Big)\vee \mathbb{H}(\alpha\circledast(\beta\rightarrow \sup\Phi))$

$=\Big(\gamma\circledast \bigvee\limits_{Q\in P(\mathbb{U})}\Phi(Q)\circledast (\mathbb{H}(\neg \mathbb{H}(\neg Q))\wedge Q)\Big)\vee \Big( \alpha\circledast \Big(\beta\rightarrow \bigvee\limits_{Q\in P(\mathbb{U})}\Phi(Q)\circledast (\mathbb{H}(Q)\vee Q\vee \mathbb{H}\mathbb{H}(Q))\Big)\Big)$.

(3) There exists an unique {\it Lvr} $\mathbb{R}$ on $\mathbb{U}$ s.t. $\mathbb{L}=\underline{\mathbb{R}}$ iff $\mathbb{L}$ fulfills the single axiom:  $\forall \Phi: P(\mathbb{U})\longrightarrow L$, $\forall \alpha\in L$,

$\mathbb{L}(\mathbb{U}\wedge (\alpha\rightarrow \inf\Phi))=\bigwedge\limits_{Q\in P(\mathbb{U})}\mathbb{U}\wedge\Big((\alpha\circledast \Phi(Q))\rightarrow \mathbb{L}(Q)\Big).$

(4) There exists an uniquely equivalent {\it Lvr} $\mathbb{R}$ on $\mathbb{U}$ s.t. $\mathbb{L}=\underline{\mathbb{R}}$ iff $\mathbb{L}$ fulfills the single axiom:   $\forall \Phi: P(\mathbb{U})\longrightarrow L$, $\forall \alpha, \beta\in L$,

$\mathbb{L}\Big(\bigwedge\limits_{Q\in P(\mathbb{U})}\mathbb{U}\wedge (\beta\circledast \Phi(Q))\rightarrow \neg \mathbb{L}(\neg Q)\Big)\wedge \mathbb{L}\Big(\mathbb{U}\wedge(\alpha\rightarrow \inf\Phi)\Big)$

$=\bigwedge\limits_{Q\in P(\mathbb{U})}\mathbb{U}\wedge \Big((\beta\circledast\Phi(Q))\rightarrow (\mathbb{L}(\neg \mathbb{L}(\neg Q))\vee Q)\Big)\wedge \Big( (\alpha\circledast \Phi(Q))\rightarrow (\mathbb{L}(Q)\wedge Q\wedge \mathbb{L}\mathbb{L}(Q))\Big)$.
\end{proposition}

\begin{remark} The axiomatic approaches 1 and 2 presented in this subsection are an extension of the axiomatic approaches 1 and 2 for $L$-rough set in Subsection 2.1. To date, an axiomatic approach 3 for $L$-valued rough set has not been provided. We will explore this topic further in the sections that follow. \end{remark}




\section{A novel axiomatic approach to {\it LVHRAO}s}

This section provides the notion of inner product (i.e., the degree of intersection) of $L$-subsets in $\mathbb{U}$, and applies it in novel  axiomatizations  on {\it LVHRAO} and  some special {\it LVHRAO}s associated with reflexive, symmetric, transitive, Euclidean and mediate {\it Lvr}s and their compositions.

\subsection{The inner product of $L$-subsets in $\mathbb{U}$ and the upper inverse mapping of $\mathbb{H}: P(\mathbb{U})\longrightarrow P(\mathbb{U})$}


\begin{definition}
Define a mapping  $\mathbb{I}: P(\mathbb{U})\times P(\mathbb{U})\longrightarrow L$  by $\forall M, Q\in P(\mathbb{U})$:
$$\mathbb{I}(M,Q)=\bigvee\limits_{d\in X}\Big(M(d)\circledast \big(\mathbb{U}(d)\rightarrow  Q(d)\big)\Big).$$  Then  $ \mathbb{I}(M,Q)$ is named the  inner product (i.e., the degree of intersection) of $M,Q$.
\end{definition}

\begin{definition} The mapping  $\mathbb{H}^{-1}$ associated with $\mathbb{H}: P(\mathbb{U})\longrightarrow   P(\mathbb{U})$ is called the upper inverse mapping of $\mathbb{H}$, where $\forall Q\in  P(\mathbb{U}), d\in X$: $ \mathbb{H}^{-1}(Q)(d)=\mathbb{I}\Big(\mathbb{U}(d)\wedge \mathbb{H}(\mathbb{U}_{\{d\}}), Q\Big).$
\end{definition}

\begin{remark} Let $\mathbb{U}=1$ and $M,Q\in P(\mathbb{U})$.

(1) Note that $$\mathbb{I}(M,Q)=\bigvee\limits_{d\in X}\Big(M(d)\circledast \big(\mathbb{U}(d)\rightarrow  Q(d)\big)\Big)=\bigvee\limits_{d\in X}\Big(M(d)\circledast Q(d)\Big).$$Hence $\mathbb{I}(M,Q)$ degenerates into the inner product (i.e., the degree of intersection) in  \cite{Bao18}.

(2) Note that
$$\mathbb{H}^{-1}(Q)(d)=\mathbb{I}\Big(\mathbb{U}(d)\wedge \mathbb{H}(\mathbb{U}_{\{d\}}), Q\Big)=\bigvee\limits_{d\in X}\Big(\mathbb{H}(\chi_{\{d\}})(d)\circledast Q(d)\Big).$$
Hence $\mathbb{H}^{-1}$  degenerates into the upper inverse mapping in \cite{Bao18}.
\end{remark}

\begin{example} Let $X=\{a,b,c,d\}$, $\mathbb{U}=\frac{0.2}{a}+\frac{0.7}{b}+\frac{0.3}{c}+\frac{0.8}{d}$ and $\mathbb{H}: P(\mathbb{U})\longrightarrow P(\mathbb{U})$ be defined by $\forall Q\in P(\mathbb{U}), \mathbb{H}(Q)=Q$.

 (1) Consider $L=[0,1]$ and $\circledast=\wedge$, then $([0,1],\wedge)$ forms a GL-quantale. $\forall \alpha, \beta\in [0,1]$, known that $\alpha\rightarrow \beta=\left\{\begin{array}{ll}1,& \alpha\leq \beta,\\ \beta,&otherwise.\end{array}\right.$.

Take $M=\frac{0.2}{a}+\frac{0.5}{b}+\frac{0.3}{c}+\frac{0.6}{d}$ and $Q=\frac{0.2}{a}+\frac{0.5}{b}+\frac{0.3}{c}+\frac{0.5}{d}$, then $M, Q\in P(\mathbb{U})$. So, $$\mathbb{I}(M,Q)=0.2\vee 0.5\vee 0.3 \vee 0.5=0.5, \mathbb{H}^{-1}(Q)=\frac{0.2}{a}+\frac{0.5}{b}+\frac{0.3}{c}+\frac{0.5}{d}.$$

(2) Consider $L=[0,1]$ and $\circledast$ being defined through $\forall \alpha, \beta\in [0,1]$, $\alpha\circledast \beta=\max\{\alpha+\beta-1, 0\}$,  then $([0,1],\circledast)$ forms a GL-quantale, and $\alpha\rightarrow \beta=\min\{1,1-\alpha+\beta\}$.

Take $M=\frac{0.2}{a}+\frac{0.5}{b}+\frac{0.3}{c}+\frac{0.6}{d}$ and $Q=\frac{0.2}{a}+\frac{0.5}{b}+\frac{0.3}{c}+\frac{0.5}{d}$, then $M, Q\in P(\mathbb{U})$. So, $$\mathbb{I}(M,Q)=0.2\vee 0.3\vee 0.3 \vee 0.3=0.3,\mathbb{H}^{-1}(Q)=\frac{0.2}{a}+\frac{0.5}{b}+\frac{0.3}{c}+\frac{0.5}{d}.$$

\end{example}

The subsequent  propositions give the frequently used properties of the proposed two notions.

\begin{proposition} \label{lemmaOI} Let $W, V, V_i(i\in I)\in P(\mathbb{U})$.

{\rm (1)} $\forall D\in P(\mathbb{U}), \mathbb{I}\Big(D, W\Big)\leq \mathbb{I}\Big(D, V\Big)$ implies $W\leq V$, furthermore,  $\mathbb{I}\Big(D, W\Big)= \mathbb{I}\Big(D, V\Big)$ implies $W= V$.

{\rm (2)} $\mathbb{I}\Big(V, W\Big)=\mathbb{I}\Big(W, V\Big)$.

{\rm (3)} $\forall \alpha, \beta \in  L$ with $\beta\leq \alpha$ and $\bigvee\limits_{u\in X}V(u)\leq \alpha$, then
$\mathbb{I}\Big(W,  \beta\circledast (\alpha\rightarrow V)\Big)=\beta\circledast\Big(\alpha\rightarrow \mathbb{I}\Big(W, V\Big)\Big).$

{\rm (4)} $\mathbb{I}\Big(W,  \bigvee\limits_{i} V_i\Big)=\bigvee\limits_{i}\mathbb{I}\Big(W, V_i\Big)$.
\end{proposition}

\begin{proof}  (1) $\forall d\in X$,

\begin{eqnarray*}\mathbb{I}\Big(\mathbb{U}_{\{d\}}, W\Big)=\bigvee_{h\in X}\Big(\mathbb{U}_{\{d\}}(h)\circledast \big(\mathbb{U}(h)\rightarrow W(h)\big)\Big)=\mathbb{U}(d)\circledast \big(\mathbb{U}(d)\rightarrow W(d)\big)=W(d).\end{eqnarray*}Similarly, $\mathbb{I}\Big(\mathbb{U}_{\{d\}}, V\Big)=V(d)$. So $\mathbb{I}\Big(\mathbb{U}_{\{d\}}, W\Big)\leq \mathbb{I}\Big(\mathbb{U}_{\{d\}}, V\Big)$ implies $W\leq V$. The further part follows similar.

(2) It is conclude by Lemma \ref{jb}(2).

(3) From  $\bigvee\limits_{u\in X}V(u)\leq \alpha$ and $\beta\leq \alpha$ we realize  $\beta\circledast (\alpha\rightarrow V)= V\circledast (\alpha\rightarrow \beta)\in  P(\mathbb{U})$ and  $\mathbb{I}\Big(W,  V\Big)\leq \alpha$. Then
\begin{eqnarray*}\mathbb{I}\Big(W,  \beta\circledast (\alpha\rightarrow V)\Big)&=&
 \mathbb{I}\Big(W, V\circledast (\alpha\rightarrow \beta)\Big)\\
&=&(\alpha\rightarrow \beta)\circledast \mathbb{I}\Big(W,  V\Big)=\beta\circledast\Big(\alpha\rightarrow \mathbb{I}\Big(W, V\Big)\Big).\end{eqnarray*}

(4) It holds by
\begin{eqnarray*}\mathbb{I}\Big(W,  \bigvee_{i} V_i\Big)&=& \mathbb{I}\Big(\bigvee_{i} V_i, W\Big)=\bigvee\limits_{d\in X}\Big(\bigvee_{i} V_i(d)\circledast \big(\mathbb{U}(d)\rightarrow  W(d)\big)\Big)=\bigvee_{i}\mathbb{I}\Big(W, V_i\Big).\qedhere\end{eqnarray*}
 \end{proof}

Let  $\mathbb{H}: P(\mathbb{U})\longrightarrow   P(\mathbb{U})$ be  a mapping. We introduce the following the following notations.

(H0) $ \mathbb{H}(\mathbb{U}_{\{d\}})\leq \mathbb{U}(d)$ for any $d\in X$.

(H1)  $\mathbb{H}   (\beta \circledast (\alpha\rightarrow Q))=\beta \circledast(\alpha\rightarrow \mathbb{H} (Q))$ for any $\beta ,\alpha\in L, Q \in P(\mathbb{U})$ with  $\beta \leq \alpha$ and $\bigvee\limits_{d\in X}Q(d)\leq \alpha$.

(H2)   $\mathbb{H}   (\bigvee\limits_{i}W_i)=\bigvee\limits_{i} \mathbb{H}   (W_i)$ for any $W_i(i\in I)\subseteq P(\mathbb{U})$.

\begin{remark} \label{remger} From Proposition \ref{theoremL1} (1) we know that $\mathbb{H}$ is an {\it LVHRAO} iff it fulfills (H1) and (H2).\end{remark}

\begin{proposition} \label{lemmafn} Suppose  $\mathbb{H}: P(\mathbb{U})\longrightarrow   P(\mathbb{U})$ is a mapping.

 {\rm (1)}  $\forall h,d \in X$, $\mathbb{H}^{-1} (\mathbb{U}_{\{d\}})(h)=\mathbb{U}(h)\wedge \mathbb{H} (\mathbb{U}_{\{h\}})(d).$

 {\rm (2)}  If $\mathbb{H}$ satisfies (H1), then it  fulfills (H0).

 {\rm (3)}  If $\mathbb{H}$ satisfies (H0), then  $\mathbb{H}^{-1}(M)(h)=\mathbb{I}\Big(\mathbb{H}(\mathbb{U}_{\{h\}}), M\Big)$ for any $M\in P(\mathbb{U}), h\in X$.

  {\rm (4)}  If $\mathbb{H}$ satisfies (H0), then  $\mathbb{H}^{-1} (\mathbb{U}_{\{d\}})(h)=\mathbb{H} (\mathbb{U}_{\{h\}})(d)$ for any $h,d \in X$.
\end{proposition}

\begin{proof} (1) It holds by
\begin{eqnarray*} \mathbb{H}^{-1}\Big(\mathbb{U}_{\{d\}}\Big)(h)&=&\bigvee_{g\in X} \Big( \mathbb{U}_{\{d\}}(g)\circledast \big(\mathbb{U}(g)\rightarrow \mathbb{U}(h)\wedge \mathbb{H}(\mathbb{U}_{\{h\}})(g)\big)\Big)\\
&=& \mathbb{U}(d)\circledast \big(\mathbb{U}(d)\rightarrow \mathbb{U}(h)\wedge \mathbb{H}(\mathbb{U}_{\{h\}})(d)\big)=\mathbb{U}(h)\wedge \mathbb{H}(\mathbb{U}_{\{h\}})(d). \end{eqnarray*}

(2) $\forall d\in X$,  known that $\bigvee\limits_{x\in X}\mathbb{U}_{\{d\}}(x)\leq \mathbb{U}(d)$. Then it follows by (H1) that
\begin{eqnarray*}& & \mathbb{H}\Big(\mathbb{U}(d)\circledast (\mathbb{U}(d)\rightarrow \mathbb{U}_{\{d\}})\Big)=\mathbb{U}(d)\circledast\Big(\mathbb{U}(d)\rightarrow \mathbb{H} (\mathbb{U}_{\{d\}})\Big) \\
&\Longleftrightarrow&\mathbb{H}( \mathbb{U}_{\{d\}})=\mathbb{U}(d)\wedge \mathbb{H} (\mathbb{U}_{\{d\}}), \end{eqnarray*} which means $ \mathbb{H}(\mathbb{U}_{\{d\}})\leq \mathbb{U}(d)$, so (H0) holds.

(3) It follows by  that $\mathbb{H}^{-1}(M)(h)=\mathbb{I}\Big(\mathbb{U}(h)\wedge \mathbb{H}(\mathbb{U}_{\{h\}}), M\Big)=\mathbb{I}\Big(\mathbb{H}(\mathbb{U}_{\{h\}}), M\Big)$.

(4) It holds by (1). \end{proof}

\subsection{The novel single axiom characterizations on {\it LVHRAO}s}

This subsection provides the single axiom characterizations on {\it LVHRAO} and  some special {\it LVHRAO}s associated with reflexive, symmetric, transitive, Euclidean and mediate {\it Lvr}s.

Unless otherwise specified, we always assume that $\mathbb{H}: P(\mathbb{U})\longrightarrow   P(\mathbb{U})$ to be a mapping.

\begin{lemma} \label{lemmafen}  For any $Q\in  P(\mathbb{U})$, $Q=\bigvee\limits_{b\in X}\Big(Q(b)\circledast \big(\mathbb{U}(b)\rightarrow \mathbb{U}_{\{b\}}\big)\Big)$.
\end{lemma}

\begin{proof} It is concluded from that $\forall l\in X$, \begin{eqnarray*}\bigvee\limits_{b\in X}\Big(Q(b)\circledast \big(\mathbb{U}(b)\rightarrow \mathbb{U}_{\{b\}}\big)\Big)(l)&=&\bigvee\limits_{b\in X}\Big(Q(b)\circledast \big(\mathbb{U}(b)\rightarrow \mathbb{U}_{\{b\}}(l)\big)\Big)\\
&=&\bigvee\limits_{b\in X}\Big(\mathbb{U}_{\{b\}}(l)\circledast \big(\mathbb{U}(b)\rightarrow Q(b)\big)\Big)\\
&=&\mathbb{U}(l)\circledast \big(\mathbb{U}(l)\rightarrow Q(l)\big)=Q(l).\qedhere\end{eqnarray*}\end{proof}

\begin{theorem} \label{sutheorembg}  $\mathbb{H}$ is an {\it LVHRAO} iff it fulfills (H).

(H) $ \forall Q,  M \in P(\mathbb{U})$,
$ \mathbb{I}\Big(Q, \mathbb{H} (M)\Big)=\mathbb{I}\Big(M,  \mathbb{H}^ {-1}(Q)\Big).$
\end{theorem}

\begin{proof} From Remark \ref{remger} we only to examine   (H1)+(H2) $\Longleftrightarrow$ (H).

(1) (H)$\Longrightarrow$ (H1), (H2).

(i) Let $\alpha, \beta \in  L, Q,  M \in  P(\mathbb{U})$ with $\alpha\leq \beta $ and $\bigvee\limits_{b\in X}Q(b)\leq \beta $.

For any $d\in X$,  $\mathbb{H}(Q)(d)= \mathbb{I}\Big(\mathbb{U}_{\{d\}}, \mathbb{H} (Q)\Big)\stackrel{\rm (H)}{=}\mathbb{I}\Big(Q, \mathbb{H}^{-1}(\mathbb{U}_{\{d\}})\Big)\leq  \beta ,$ so $\bigvee \mathbb{H}(Q)\leq \beta $. Then
\begin{eqnarray*}\mathbb{I}\Big(M,  \mathbb{H} \big(\alpha\circledast (\beta \rightarrow Q)\big)\Big)&\stackrel{\rm (H)}{=}& \mathbb{I}\Big(  \alpha\circledast (\beta \rightarrow Q),   \mathbb{H}^ {-1}(M)\Big)\ {\rm by \ Proposition\   \ref{lemmaOI} (3)}
\\
&=&\alpha\circledast\Big(\beta \rightarrow \mathbb{I}\Big(Q, \mathbb{H}^ {-1}(M)\Big)\Big)\\
&\stackrel{\rm (H)}{=}&\alpha\circledast\Big(\beta \rightarrow \mathbb{I}\Big(M, \mathbb{H} (Q)\Big)\\
&=& \mathbb{I}\Big(M,\alpha\circledast(\beta \rightarrow \mathbb{H} (Q))\Big).\end{eqnarray*}From Proposition \ref{lemmaOI}(1), we obtain $ \mathbb{H} (\alpha\circledast(\beta \rightarrow Q))= \alpha\circledast(\beta \rightarrow \mathbb{H} (Q))$. Hence (H1) holds.

(ii) $\forall  M, Q_i(i\in I)\in  P(\mathbb{U})$, by (H) and  Proposition \ref{lemmaOI} (4)
\begin{eqnarray*}\mathbb{I}\Big(    M   ,  \mathbb{H} (\bigvee_{i} Q_i)\Big)&=& \mathbb{I}\Big(  \bigvee_{i} Q_i,   \mathbb{H}^ {-1}(  M   )\Big)=\bigvee_{i}\mathbb{I}\Big(M,  \mathbb{H} (Q_i)\Big)=\mathbb{I}\Big(    M   , \bigvee_{i} \mathbb{H} (Q_i)\Big).\end{eqnarray*}This shows   $ \mathbb{H} (\bigvee\limits_{i} Q_i)= \bigvee\limits_{i}  \mathbb{H} ( Q_i)$, i.e., (H2) holds.

(2) (H1), (H2) $\Longrightarrow$ (H). $\forall Q, M   \in P(\mathbb{U})$,  \begin{eqnarray*}\mathbb{I}\Big(M,  \mathbb{H}^ {-1}(Q)\Big)&=&\mathbb{I}\Big( \mathbb{H}^ {-1}(Q), M\Big)\\
&=&\bigvee_{z\in X} \Big(\mathbb{H}^ {-1}(Q)(z)\circledast\big(\mathbb{U}(z)\rightarrow M(z)\big)\Big)   \ {\rm by \ Proposition\   \ref{lemmafn}}(2),(3)\\
&=&\bigvee_{z\in X} \Big(\bigvee_{b\in X}\Big[\mathbb{H}(\mathbb{U}_{\{z\}})(b)\circledast \big(\mathbb{U}(b)\rightarrow Q(b)\big)\Big]\circledast\big(\mathbb{U}(z)\rightarrow M(z)\big)\Big)\\
&=&\bigvee_{b\in X} \Big(\bigvee_{z\in X}\Big[ \mathbb{H}(\mathbb{U}_{\{z\}})(b)\circledast \big(\mathbb{U}(z)\rightarrow M(z)\big)\Big]\circledast\big(\mathbb{U}(b)\rightarrow Q(b)\big)\Big)\\
&=&\bigvee_{b\in X} \Big(\bigvee_{z\in X}\Big[M(z)\circledast \big(\mathbb{U}(z)\rightarrow  \mathbb{H}(\mathbb{U}_{\{z\}})(b)\big)\Big]\circledast\big(\mathbb{U}(b)\rightarrow Q(b)\big)\Big) \ {\rm by\ (H1), (H2), Lemma\   \ref{lemmafen}}\\
&=&\bigvee_{b\in X} \Big(\mathbb{H}(M)(b)\circledast\big(\mathbb{U}(b)\rightarrow Q(b)\big)\Big)\\
&=&\bigvee_{b\in X} \Big(Q(b)\circledast\big(\mathbb{U}(b)\rightarrow \mathbb{H}(M)(b)\big)\Big)=\mathbb{I}\Big(Q, \mathbb{H} (M)\Big).  \qedhere\end{eqnarray*}
\end{proof}

To facilitate the expression, we give the following concepts.

\begin{definition} A mapping  $\mathbb{H}: P(\mathbb{U})\longrightarrow P(\mathbb{U})$ is named a reflexive {\it LVHRAO}  if there exists a reflexive {\it Lvr} $\mathbb{R}$ on $\mathbb{U}$ s.t. $\mathbb{H}=\overline{\mathbb{R}}$. Similarly, one can define many other special  {\it LVHRAO}.

\end{definition}

\begin{remark} Next, we will describe some special {\it LVHRAO}s by strengthening axiom (H). For readers to quickly understand the operators described by these axioms, we will abbreviate ``reflexive, transitive, symmetric, Euclidean and mediate" as ``R, T, S, E and M" respectively. Then, readers can understand that axiom (HR)  describes reflexive {\it LVHRAO}, (HRT) describes reflexive and transitive {\it LVHRAO}, (HRTS) describes reflexive, transitive and symmetric {\it LVHRAO}.
\end{remark}

\begin{proposition} \label{opropr}  $\mathbb{H}$ is a reflexive {\it LVHRAO} iff it fulfills (H1), (H2) and (H3).

(H3) $\forall W\in  P(\mathbb{U})$, $\mathbb{H}^{-1}(W)\geq W$.
\end{proposition}

\begin{proof}
  $\Longrightarrow$. Let $\mathbb{H} =\overline{\mathbb{R}}$ for a reflexive {\it Lvr} $\mathbb{R}$ on $\mathbb{U}$. We only examine (H3). Indeed, $\forall W\in  P(\mathbb{U})$, $b\in X$, by Proposition \ref{lemmafn}(2),(3) \begin{eqnarray*}\mathbb{H}^{-1}(W)(b)&=&\bigvee\limits_{z\in X}\Big(\overline{\mathbb{R}}(\mathbb{U}_{\{b\}})(z)\circledast (\mathbb{U}(z)\rightarrow W(z)) \Big)\\
 &=&\bigvee\limits_{z\in X}\Big(\mathbb{R}(b,z)\circledast (\mathbb{U}(z)\rightarrow W(z)) \Big)\\
 &\geq&\mathbb{R}(b,b)\circledast (\mathbb{U}(b)\rightarrow W(b)) \ {\rm by \ reflexivity}\\
 &\geq& \mathbb{U}(b)\circledast (\mathbb{U}(b)\rightarrow W(b))= W(b).\end{eqnarray*}

$\Longleftarrow$. By Remark \ref{remger}, $\mathbb{H} =\overline{\mathbb{R}}$ for an {\it Lvr} $\mathbb{R}$. Then $\forall b\in X$,  \begin{eqnarray*}\mathbb{U}(b)&=&\mathbb{U}_{\{b\}}(b)\stackrel{(H3)}{\leq} \mathbb{H}^{-1}(\mathbb{U}_{\{b\}})(b)   \ {\rm by \ Proposition\   \ref{lemmafn}}(1)\\
&=&\mathbb{U}(b)\wedge \overline{\mathbb{R}}(\mathbb{U}_{\{b\}})(b)=\mathbb{R}(b,b).\end{eqnarray*} So, $\mathbb{R}$ is reflexive.
\end{proof}

\begin{theorem} \label{osutheoremr} $\mathbb{H}$ is a reflexive {\it {\it LVHRAO}} iff it fulfills (HR).

(HR) $\forall M, Q\in P(\mathbb{U})$, $\mathbb{I}\Big(M,\mathbb{H}(Q)\Big)=\mathbb{I}\Big(Q, M\vee \mathbb{H}^{-1}(M)\Big).$
\end{theorem}

\begin{proof}  By Theorem \ref{sutheorembg} and Proposition \ref{opropr}, we only need examine (H)+(H3) $\Longleftrightarrow$ (HR).

(1) (H)+(H3) $\Longrightarrow$ (HR). Indeed  $$\mathbb{I}\Big(Q, M\vee \mathbb{H}^{-1}(M)\Big)\stackrel{\rm (H3)}{=}\mathbb{I}\Big(Q, \mathbb{H}^{-1}(M)\Big)\stackrel{\rm (H)}{=}\mathbb{I}\Big(M,\mathbb{H}(Q)\Big),$$ i.e., (HR) holds.

(2) (HR) $\Longrightarrow$ (H)+(H3).

(i) (HR) $\Longrightarrow$ (H0). $\forall d, b\in X$,  put $M=\mathbb{U}_{\{b\}}, Q=\mathbb{U}_{\{d\}}$ in (HR) then
\begin{eqnarray*}& &\mathbb{I}\Big(\mathbb{U}_{\{b\}},\mathbb{H}(\mathbb{U}_{\{d\}})\Big)\stackrel{\rm (HR)}{=}\mathbb{I}\Big(\mathbb{U}_{\{d\}}, \mathbb{U}_{\{b\}}\vee \mathbb{H}^{-1}(\mathbb{U}_{\{b\}})\Big)\\
&\Longleftrightarrow& \mathbb{H}(\mathbb{U}_{\{d\}})(b)=\mathbb{U}(d)\wedge \Big(\mathbb{U}_{\{b\}}(d)\vee \mathbb{H}^{-1}(\mathbb{U}_{\{b\}})(d)\Big),\end{eqnarray*}which means $\mathbb{H}(\mathbb{U}_{\{d\}})(b)\leq \mathbb{U}(d)$, i.e., (H0) holds.

(ii) (HR) $\Longrightarrow$ (H3). $\forall u\in X$,  put $Q=\mathbb{U}_{\{u\}}$ in (HR) then

$$\mathbb{H}^{-1}(M)(u)\stackrel{\rm Proposition\   \ref{lemmafn}(3)}{=}\mathbb{I}\Big(M,\mathbb{H}(\mathbb{U}_{\{u\}})\Big)\stackrel{\rm (HR)}{=}\mathbb{I}\Big(\mathbb{U}_{\{u\}}, M\vee \mathbb{H}^{-1}(M)\Big)\geq  M(u),$$ so (H3) gains.

(iii) Using (H3) to (HR) one get (H).
\end{proof}

\begin{proposition} \label{opropt}  $\mathbb{H}$ is a transitive LVHRAO iff it fulfills (H1), (H2) and (H4).

 (H4) $\forall M \in P(\mathbb{U})$, $\mathbb{H}^{-1}\mathbb{H}^{-1}(M)\leq \mathbb{H}^{-1} (M)$.
\end{proposition}

\begin{proof} $\Longrightarrow$.  Let $\mathbb{H} =\overline{\mathbb{R}}$ for a transitive {\it Lvr} $\mathbb{R}$ on $\mathbb{U}$. We only examine (H4). Indeed, $\forall M\in  P(\mathbb{U})$, $d\in X$,
 \begin{eqnarray*}\mathbb{H}^{-1}\mathbb{H}^{-1}(M)(d)&=&\overline{\mathbb{R}} ^{-1}\overline{\mathbb{R}}^{-1}(M)(d)   \ {\rm by \ Proposition\   \ref{lemmafn}(2),(3)}\\
&=&\bigvee_{z\in X}\Big(\bigvee_{b\in X}\Big[\overline{\mathbb{R}}(\mathbb{U}_{\{z\}})(b)\circledast(\mathbb{U}(b)\rightarrow M(b))\Big]\circledast \Big[\mathbb{U}(z)\rightarrow \mathbb{R}(d,z)\Big]\Big)\\
&=&\bigvee_{z\in X}\bigvee_{b\in X}\Big((\mathbb{U}(b)\rightarrow M(b))\circledast\mathbb{R} (z,b)\circledast \Big[\mathbb{U}(z)\rightarrow \mathbb{R}(d,z)\Big]\Big)\\
&=&\bigvee_{z\in X}\bigvee_{b\in X}\Big((\mathbb{U}(b)\rightarrow M(b))\circledast\mathbb{R} (d,z)\circledast \Big[\mathbb{U}(z)\rightarrow \mathbb{R}(z,b)\Big]\Big) \ {\rm by \ transitive}\\
&\leq&\bigvee_{b\in X}\Big((\mathbb{U}(b)\rightarrow M(b))\circledast \mathbb{R}(d,b)\Big)\\
&=&\bigvee_{b\in X}\Big((\mathbb{U}(b)\rightarrow M(b))\circledast\overline{\mathbb{R}}(\mathbb{U}_{\{d\}})(b)\Big)   \ {\rm by \ Proposition\   \ref{lemmafn}(2),(3)}\\
&=&\mathbb{H}^{-1}(M)(d).\end{eqnarray*}

$\Longleftarrow$. By Remark \ref{remger}, $\mathbb{H} =\overline{\mathbb{R}}$ for an {\it Lvr} $\mathbb{R}$. Then  $\forall b, z\in X$,  \begin{eqnarray*}\mathbb{R}(z,b)&=& \overline{\mathbb{R}}^{-1}(\mathbb{U}_{\{b\}})(z)\\
&\stackrel{\rm (H4)}{\geq}& \overline{\mathbb{R}}^{-1}\overline{\mathbb{R}} ^{-1}(\mathbb{U}_{\{b\}})(z)   \ {\rm by \ Proposition\   \ref{lemmafn}(2),(3),(4)}\\
&=&\bigvee_{d\in X}\Big(\overline{\mathbb{R}}(\mathbb{U}_{\{z\}})(d)\circledast \Big[\mathbb{U}(d)\rightarrow \overline{\mathbb{R}} ^{-1} (\mathbb{U}_{\{b\}})(d)\Big]\Big)\\
&=&\bigvee_{d\in X}\Big(\mathbb{R}(z,d)\circledast \Big[\mathbb{U}(d)\rightarrow \mathbb{R}(d,b)\Big]\Big).\end{eqnarray*} So, $\mathbb{R}$ is transitive.
\end{proof}

\begin{theorem} \label{osutheoremt}  $\mathbb{H}$ is  a transitive  {\it LVHRAO} iff it satisfies (HT).

(HT) $\forall M, Q\in P(\mathbb{U})$, $ \mathbb{I}\Big(M,\mathbb{H}(Q)\Big)=\mathbb{I}\Big(Q, \mathbb{H}^{-1}\mathbb{H}^{-1}(M)\vee \mathbb{H}^{-1}(M)\Big).$
\end{theorem}

\begin{proof} From Theorem \ref{sutheorembg} and Proposition \ref{opropt} we only need examine (H)+(H4) $\Longleftrightarrow$ (HT).

(1) (H)+(H4) $\Longrightarrow$ (HT). Indeed $$\mathbb{I}\Big(Q, \mathbb{H}^{-1}(M)\vee \mathbb{H}^{-1}\mathbb{H}^{-1}(M)\Big)\stackrel{\rm (H4)}{=}\mathbb{I}\Big(Q, \mathbb{H}^{-1}(M)\Big)\stackrel{\rm (H)}{=}\mathbb{I}\Big(M,\mathbb{H}(Q)\Big),$$ i.e., (HT) holds.

(1) (HT) $\Longrightarrow$ (H)+(H4).

Analogizing the verification of (i) in Theorem \ref{osutheoremr}, we get that (HT) $\Longrightarrow$ (H0). Then $\forall b\in X$,  take $Q=\mathbb{U}_{\{b\}}$ in (HT) one get $$\mathbb{H}^{-1}(M)(b)\stackrel{\rm Proposition\   \ref{lemmafn}(3)}{=}\mathbb{I}\Big(M,\mathbb{H}(\mathbb{U}_{\{b\}})\Big)\stackrel{\rm (HT)}{=}\mathbb{I}\Big(\mathbb{U}_{\{b\}}, \mathbb{H}^{-1}\mathbb{H}^{-1}(M)\vee \mathbb{H}^{-1}(M)\Big)\geq  \mathbb{H}^{-1}\mathbb{H}^{-1}(M)(b),$$ so (H4) gains. Using (H4) to (HT) one get (H).
\end{proof}

\begin{proposition} \label{oprops}    $\mathbb{H}$ is a symmetric LVHRAO iff it fulfills (H1), (H2) and (H5).

(H5) $\forall M\in P(\mathbb{U})$, $\mathbb{H}^{-1}(M)=\mathbb{H}(M)$.
\end{proposition}

\begin{proof}
$\Longrightarrow$. Let $\mathbb{H} =\overline{\mathbb{R}}$ for a symmetric {\it Lvr} $\mathbb{R}$ on $\mathbb{U}$. We only examine (H5). $\forall M\in  P(\mathbb{U})$, $z\in X$,
 \begin{eqnarray*}\mathbb{H}^{-1}(M)(z)&=&\bigvee_{b\in X}\Big(\overline{\mathbb{R}}(\mathbb{U}_{\{z\}})(b)\circledast \Big[\mathbb{U}(b)\rightarrow M(b)\Big]\Big)\\
&=&\bigvee_{b\in X}\Big(\mathbb{R}(z,b)\circledast \Big[\mathbb{U}(b)\rightarrow M(b)\Big]\Big) \ {\rm by \ symmetry}\\
&=&\bigvee_{b\in X}\Big(\mathbb{R}(b,z)\circledast \Big[\mathbb{U}(b)\rightarrow M(b)\Big]\Big)=\mathbb{H}(M)(z).\end{eqnarray*}

$\Longleftarrow$. By Remark \ref{remger}, $\mathbb{H} =\overline{\mathbb{R}}$ for an {\it Lvr} $\mathbb{R}$. Then $\forall z,b\in X$,  \begin{eqnarray*}\mathbb{R}(b,z)=\overline{\mathbb{R}}(\mathbb{U}_{\{b\}})(z)\stackrel{\rm (H5)}{=}\overline{\mathbb{R}}^{-1}(\mathbb{U}_{\{b\}})(z)\stackrel{\rm Proposition \ref{lemmafn}(2),(4)}{=}\overline{\mathbb{R}}(\mathbb{U}_{\{z\}})(b)=\mathbb{R}(z,b).\end{eqnarray*} So, $\mathbb{R}$ is symmetric.
\end{proof}

\begin{theorem} \label{osutheorems}  $\mathbb{H}$ is  a symmetric   {\it LVHRAO} iff it satisfies (HS).

(HS) $\forall M, Q\in P(\mathbb{U})$, $\mathbb{I}\Big(M,\mathbb{H}(Q)\Big)=\mathbb{I}\Big(Q,  \mathbb{H}(M)\Big).$
\end{theorem}

\begin{proof} From Theorem \ref{sutheorembg} and Proposition \ref{oprops} we only need examine (H)+(H5) $\Longleftrightarrow$ (HS).

(1) (H)+(H5) $\Longrightarrow$ (HS).  Indeed $$\mathbb{I}\Big(Q, \mathbb{H}(M)\Big)\stackrel{\rm (H5)}{=}\mathbb{I}\Big(Q, \mathbb{H}^{-1}(M)\Big)\stackrel{\rm (H)}{=}\mathbb{I}\Big(M,\mathbb{H}(Q)\Big).$$

(2) (HS) $\Longrightarrow$ (H)+(H5). Analogizing the verification  of (i) in Theorem \ref{osutheoremr}, we have that (HT) $\Longrightarrow$ (H0). Then $\forall b\in X$,  take $Q=\mathbb{U}_{\{b\}}$ in (HS) one get $$\mathbb{H}^{-1}(M)(b)\stackrel{\rm Proposition\   \ref{lemmafn}(3)}{=}\mathbb{I}\Big(M,\mathbb{H}(\mathbb{U}_{\{b\}})\Big)\stackrel{\rm (HS)}{=}\mathbb{I}\Big(\mathbb{U}_{\{b\}}, \mathbb{H}(M)\Big)= \mathbb{H}(M)(b),$$  so  (H5) gains. Using (H5) to (HS) one get (H).
\end{proof}

\begin{definition} \cite{XW22}  An {\it Lvr} $\mathbb{R}$ is referred Euclidean whenever $\forall a,d,h\in X$, $\mathbb{R}(d,a)\circledast (\mathbb{U}(d)\rightarrow \mathbb{R}(d,h))\leq \mathbb{R}(a,h)$.
\end{definition}

If $\mathbb{U}=1$, then  an {\it Lvr} $\mathbb{R}$ on $\mathbb{U}$ degenerates into an $L$-relation on $X$. The Euclidean condition reduces into $\mathbb{R}(d,a)\circledast \mathbb{R}(d,h)\leq \mathbb{R}(a,h)$, i.e., the Euclidean condition for $L$-relation in \cite{AM05}.

\begin{example} Let $X=\{a,b,c\}$, $L=[0,1]$, $\circledast=\wedge$ and  $\mathbb{U}=\frac{0.5}{a}+\frac{0.7}{b}+\frac{0.4}{c}$. Then the mapping $\mathbb{R}:X\times X\longrightarrow [0,1]$ define below is an Euclidean {\it Lvr} on $\mathbb{U}$, where

\ \

 \begin{tabular}{l|lll}

   $\mathbb{R}(-,-)$    &  $a$ & $b$ & $c$   \\
\hline
$a$ &  0.5 & 0.2 & 0.2  \\

$b$  & 0.2 & 0.7 & 0.1    \\

$c$  & 0.2   & 0.1 & 0.4 \\
 \end{tabular}
\end{example}

\begin{proposition} \label{opropeu}  $\mathbb{H}$ is an  Euclidean LVHRAO iff it fulfills (H1), (H2) and (H6).

(H6) $\forall M \in P(\mathbb{U})$, $\mathbb{H}\mathbb{H}^{-1} (M)\leq \mathbb{H}^{-1} (M)$.
\end{proposition}

\begin{proof} $\Longrightarrow$. Let $\mathbb{H} =\overline{\mathbb{R}}$ for an  Euclidean {\it Lvr} $\mathbb{R}$ on $\mathbb{U}$. We only examine (H6). $\forall M\in  P(\mathbb{U})$, $d\in X$,
 \begin{eqnarray*}\mathbb{H}\mathbb{H} ^{-1}(M)(d)&=&\bigvee_{b\in X}\Big(\overline{\mathbb{R}}^{-1}(M)(b)\circledast \Big[\mathbb{U}(b)\rightarrow \mathbb{R}(b,d)\Big]\Big)   \ {\rm from \ Proposition\   \ref{lemmafn}(2),(3)}\\
&=&\bigvee_{b\in X}\Big(\bigvee_{z\in X}\Big[\overline{\mathbb{R}}(\mathbb{U}_{\{b\}})(z)\circledast(\mathbb{U}(z)\rightarrow M(z))\Big]\circledast \Big[\mathbb{U}(b)\rightarrow \mathbb{R}(b,d)\Big]\Big)\\
&=&\bigvee_{b\in X}\Big(\bigvee_{z\in X}\Big[\mathbb{R} (b,z)\circledast(\mathbb{U}(z)\rightarrow M(z))\Big]\circledast \Big[\mathbb{U}(b)\rightarrow \mathbb{R}(b,d)\Big]\Big)\\
&=&\bigvee_{b\in X}\bigvee_{z\in X}\Big((\mathbb{U}(z)\rightarrow M(z))\circledast\mathbb{R} (b,d)\circledast \Big[\mathbb{U}(b)\rightarrow \mathbb{R}(b,z)\Big]\Big) \ {\rm by \ Euclidean}\\
&\leq&\bigvee_{z\in X}\Big((\mathbb{U}(z)\rightarrow M(z))\circledast \mathbb{R}(d,z)\Big)\\
&=&\bigvee_{z\in X}\Big((\mathbb{U}(z)\rightarrow M(z))\circledast\overline{\mathbb{R}}(\mathbb{U}_{\{d\}})(z)\Big)=\mathbb{H}^{-1}(M)(d).\end{eqnarray*}

$\Longleftarrow$. By Remark \ref{remger}, $\mathbb{H} =\overline{\mathbb{R}}$ for an {\it Lvr} $\mathbb{R}$. Then $\forall b, d\in X$, by Proposition \ref{lemmafn}
\begin{eqnarray*}\mathbb{R}(d,b)&=& \overline{\mathbb{R}}^{-1}(\mathbb{U}_{\{b\}})(d)\stackrel{\rm (H6)}{\geq} \overline{\mathbb{R}}\overline{\mathbb{R}}^{-1}(\mathbb{U}_{\{b\}})(d)\\
&=&\bigvee_{g\in X}\Big(\mathbb{R}(g,d)\circledast \Big[\mathbb{U}(g)\rightarrow \mathbb{R}(g,b)\Big]\Big).\end{eqnarray*} So, $\mathbb{R}$ is Euclidean.
\end{proof}

\begin{theorem} \label{sutheoremeu}  $\mathbb{H}$ is  an Euclidean  {\it LVHRAO} iff it fulfills  (HE).

(HE) $\forall M, Q\in P(\mathbb{U})$, $\mathbb{I}\Big(M,\mathbb{H}(Q)\Big)=\mathbb{I}\Big(Q, \mathbb{H}\mathbb{H}^{-1}(M)\vee \mathbb{H}^{-1}(M)\Big).$
\end{theorem}

\begin{proof} From Theorem \ref{sutheorembg} and Proposition \ref{opropeu} we only need examine (H)+(H6) $\Longleftrightarrow$ (HE).

(1) (H)+(H6) $\Longrightarrow$ (HE).  Indeed $$\mathbb{I}\Big(Q, \mathbb{H}^{-1}(M)\vee \mathbb{H}\mathbb{H}^{-1}(M)\Big)\stackrel{\rm (H6)}{=}\mathbb{I}\Big(Q, \mathbb{H}^{-1}(M)\Big)\stackrel{\rm (H)}{=}\mathbb{I}\Big(M,\mathbb{H}(Q)\Big).$$

(2) (HE) $\Longrightarrow$ (H)+(H6). Analogizing the verification  of (i) in Theorem \ref{osutheoremr}, we have that (HT) $\Longrightarrow$ (H0). Then $\forall b\in X$,  take $Q=\mathbb{U}_{\{b\}}$ in (HE) one get $$\mathbb{H}^{-1}(M)(b)\stackrel{\rm Proposition\   \ref{lemmafn}(3)}{=}\mathbb{I}\Big(M,\mathbb{H}(\mathbb{U}_{\{b\}})\Big)\stackrel{\rm (HE)}{=}\mathbb{I}\Big(\mathbb{U}_{\{b\}}, \mathbb{H}\mathbb{H}^{-1}(M)\vee \mathbb{H}^{-1}(M)\Big)\geq  \mathbb{H}\mathbb{H}^{-1}(M)(b),$$so  (H6) gains. Using (H6) to (HE) one get (H).
\end{proof}

\begin{definition} \cite{XW22}  An {\it Lvr} $\mathbb{R}$ on $\mathbb{U}$ is referred mediate whenever $\forall a,d\in X$, $\mathbb{R}(d,a)\leq \bigvee\limits_{h\in X}\Big(\mathbb{R}(d,h)\circledast (\mathbb{U}(h)\rightarrow \mathbb{R}(h,a))\Big)$.
\end{definition}

If $\mathbb{U}=1$, then  an {\it Lvr} $\mathbb{R}$ on $\mathbb{U}$ degenerates into an $L$-relation on $X$. The mediate condition reduces into $\mathbb{R}(d,a)\leq \bigvee\limits_{h\in X}\Big(\mathbb{R}(d,h)\circledast \mathbb{R}(h,a)\Big)$, i.e., the mediate condition for $L$-relation in \cite{JIN232}.

\begin{example} Let $X=\{a,b,c\}$, $L=[0,1]$ and  $\mathbb{U}=\frac{0.5}{a}+\frac{0.7}{b}+\frac{0.4}{c}$. Then the mapping $\mathbb{R}:X\times X\longrightarrow [0,1]$ define below is a mediate {\it Lvr} on $\mathbb{U}$, where

\ \

 \begin{tabular}{l|lll}

   $\mathbb{R}(-,-)$    &  $a$ & $b$ & $c$   \\
\hline
$a$ &  0.5 & 0.2 & 0.3  \\

$b$  & 0.1 & 0.7 & 0.1    \\

$c$  & 0.2   & 0.4 & 0.4 \\
 \end{tabular}
\end{example}

\begin{proposition} \label{opropm}  $\mathbb{H}$ is a mediate LVHRAO iff it fulfills (H1), (H2) and (H7).

 (H7) $\forall W \in P(\mathbb{U})$, $\mathbb{H}^{-1}\mathbb{H}^{-1}(W)\geq \mathbb{H}^{-1} (W)$.
\end{proposition}

\begin{proof} The proof is similar to Proposition \ref{opropt}.
\end{proof}

\begin{theorem} \label{osutheoremm}  $\mathbb{H}$ is  a mediate  {\it LVHRAO} iff it satisfies (HM).

(HM) $\forall M, Q\in P(\mathbb{U})$, $ \mathbb{I}\Big(M,\mathbb{H}(Q)\Big)=\mathbb{I}\Big(Q, \mathbb{H}^{-1}\mathbb{H}^{-1}(M)\wedge \mathbb{H}^{-1}(M)\Big).$
\end{theorem}

\begin{proof} From Theorem \ref{sutheorembg} and Proposition \ref{opropm} we only need examine (H)+(H7) $\Longleftrightarrow$ (HM). The detailed proof is similar to Theorem \ref{osutheoremt}. \end{proof}
%
%
%
%

\subsection{The single axiom characterizations of composite LVHRAOs}

This subsection provides the single axiom characterizations of LVHRAOs derived from some composition about reflexive, symmetric, transitive, Euclidean $\&$ mediate {\it Lvr}s.

We first consider the combinations of two conditions.

\begin{theorem} \label{osutheoremrt}  (1) $\mathbb{H}$ is a reflexive and transitive, i.e., preordered  {\it LVHRAO} iff it fulfills (HRT).

(HRT) $\forall M, Q\in P(\mathbb{U})$, $\mathbb{I}\Big(M,\mathbb{H}(Q)\Big)=\mathbb{I}\Big(Q, M\vee \mathbb{H}^{-1}\mathbb{H}^{-1}(M)\vee \mathbb{H}^{-1}(M)\Big).$

(2) $\mathbb{H}$ is a reflexive and symmetric, i.e., tolerance  {\it LVHRAO} iff it fulfills (HRS).

(HRS) $\forall M, Q\in P(\mathbb{U})$, $\mathbb{I}\Big(M,\mathbb{H}(Q)\Big)=\mathbb{I}\Big(Q, M\vee  \mathbb{H}(M)\Big).$
\end{theorem}

\begin{proof} (1)  By Theorem \ref{osutheoremr} and Proposition \ref{opropt}, we  need examine (HR)+(H4) $\Longleftrightarrow$ (HRT).

(i) (HR)+(H4) $\Longrightarrow$ (HRT). Indeed  $$\mathbb{I}\Big(Q, M\vee \mathbb{H}^{-1}\mathbb{H}^{-1}(M)\vee \mathbb{H}^{-1}(M)\Big)\stackrel{\rm (H4)}{=}\mathbb{I}\Big(Q, M\vee \mathbb{H}^{-1}(M)\Big)\stackrel{\rm (HR)}{=}\mathbb{I}\Big(M,\mathbb{H}(Q)\Big),$$ i.e., (HRT) holds.

(ii) (HRT) $\Longrightarrow$ (HR)+(H4).

Analogizing the verification  of (i) in Theorem \ref{osutheoremr}, we have that (HRT) $\Longrightarrow$ (H0). Then $\forall o\in X$,  put $Q=\mathbb{U}_{\{o\}}$ in (HRT)

$$\mathbb{H}^{-1}(M)(o)\stackrel{\rm Proposition\   \ref{lemmafn}(3)}{=}\mathbb{I}\Big(M,\mathbb{H}(\mathbb{U}_{\{o\}})\Big)\stackrel{\rm (HRT)}{=}\mathbb{I}\Big(\mathbb{U}_{\{o\}}, M\vee  \mathbb{H}^{-1}\mathbb{H}^{-1}(M)\vee \mathbb{H}^{-1}(M)\Big)\geq  \mathbb{H}^{-1}\mathbb{H}^{-1}(M)(o),$$ so  (H4) gains. Using (H4) to (HRT) one get (HR).

(2) By Theorem \ref{osutheoremr} and Proposition \ref{oprops}, we only need examine (HR)+(H5) $\Longleftrightarrow$ (HRS). The detailed proof is similar to (1).
\end{proof}

\begin{proposition} \label{osutheoremts} (1) $\mathbb{H}$ is a reflexive and  Euclidean  {\it LVHRAO} iff it fulfills (HRE).

(HRE) $\forall Q, N\in P(\mathbb{U})$, $\mathbb{I}\Big(Q,\mathbb{H}(N)\Big)=\mathbb{I}\Big(N, Q\vee \mathbb{H}\mathbb{H}^{-1}(Q)\vee \mathbb{H}^{-1}(Q)\Big).$

(2) $\mathbb{H}$ is  a transitive and symmetric  {\it LVHRAO} iff it satisfies (HTS).

(HTS) $\forall Q, N\in P(\mathbb{U})$, $ \mathbb{I}\Big(Q,\mathbb{H}(N)\Big)=\mathbb{I}\Big(N, \mathbb{H}\mathbb{H}(Q)\vee \mathbb{H}(Q)\Big).$

(3)  $\mathbb{H}$ is  a transitive and Euclidean  {\it LVHRAO} iff it satisfies (HTE).

(HTE) $\forall Q, N\in P(\mathbb{U})$, $ \mathbb{I}\Big(Q,\mathbb{H}(N)\Big)=\mathbb{I}\Big(N, \mathbb{H}^{-1}\mathbb{H}^{-1}(Q)\vee \mathbb{H}\mathbb{H}^{-1}(Q)\vee \mathbb{H}^{-1}(Q)\Big).$

(4) $\mathbb{H}$ is  a transitive and mediate  {\it LVHRAO} iff it satisfies (HTM).

(HTM) $\forall M, N\in P(\mathbb{U})$, $ \mathbb{I}\Big(M,\mathbb{H}(N)\Big)=\mathbb{I}\Big(N, \mathbb{H}^{-1}\mathbb{H}^{-1}(M)\Big).$

(5)  $\mathbb{H}$ is  a symmetric and Euclidean  {\it LVHRAO} iff it satisfies (HTS).

(6)  $\mathbb{H}$ is  a symmetric and mediate  {\it LVHRAO} iff it satisfies (HSM).

(HSM) $\forall Q, N\in P(\mathbb{U})$, $ \mathbb{I}\Big(Q,\mathbb{H}(N)\Big)=\mathbb{I}\Big(N, \mathbb{H}\mathbb{H}(Q)\wedge \mathbb{H}(Q)\Big).$
\end{proposition}

\begin{proof}  (1) By Theorem \ref{osutheoremr} and Proposition \ref{opropeu}, we only need examine (HR)+(H6) $\Longleftrightarrow$ (HRE). The detailed proof is similar to Theorem \ref{osutheoremrt}  (1).

(2) From Theorem \ref{osutheoremt} and Proposition \ref{oprops}, we  need examine (HT)+(H5) $\Longleftrightarrow$ (HTS).

(i) (HT)+(H5) $\Longrightarrow$ (HTS). Indeed $$\mathbb{I}\Big(N, \mathbb{H}(Q)\vee \mathbb{H}\mathbb{H}(Q)\Big)\stackrel{\rm (H5)}{=}\mathbb{I}\Big(N, \mathbb{H}^{-1}(Q)\vee \mathbb{H}^{-1}\mathbb{H}^{-1}(Q)\Big)\stackrel{\rm (HT)}{=}\mathbb{I}\Big(Q,\mathbb{H}(N)\Big),$$ i.e., (HTS) holds.

(ii) (HTS) $\Longrightarrow$ (HT)+(H5).

Analogizing the verification  of (i) in Theorem \ref{osutheoremr}, we have that (HTS) $\Longrightarrow$ (H0). Then $\forall u\in X$, $D\in P(\mathbb{U})$ $$\mathbb{H}^{-1}(D)(u)\stackrel{\rm Proposition\   \ref{lemmafn}(3)}{=}\mathbb{I}\Big(D,\mathbb{H}(\mathbb{U}_{\{u\}})\Big)\stackrel{\rm (HTS)}{=}\mathbb{I}\Big(\mathbb{U}_{\{u\}}, \mathbb{H}\mathbb{H}(D)\vee \mathbb{H}(D)\Big)\geq  \mathbb{H}(D)(u),$$ similarly, one have  $ \mathbb{H}^{-1}(D)(u)\leq  \mathbb{H}(D)(u)$, so (H5) gains. Using (H5) to (HTS) one get (HT).

(3) From Theorem \ref{osutheoremt} and Proposition \ref{opropeu} we only need examine (HT)+(H6) $\Longleftrightarrow$ (HTE). The detailed proof is similar to (2).

(4) From Theorem \ref{osutheoremt} and Proposition \ref{opropm} we only need examine (HT)+(H7) $\Longleftrightarrow$ (HTM). The detailed proof is similar to (2).

(5) It follows by Proposition \ref{osutheoremts} (2) and that if $\mathbb{R}$ is symmetric, then $\mathbb{R}$ is transitive iff $\mathbb{R}$ is Euclidean.

(6) From Theorem \ref{osutheoremm} and Proposition \ref{oprops} we only need examine (HM)+(H5) $\Longleftrightarrow$ (HSM).

(i) (HM)+(H5) $\Longrightarrow$ (HSM). Indeed $$\mathbb{I}\Big(N, \mathbb{H}(Q)\wedge \mathbb{H}\mathbb{H}(Q)\Big)\stackrel{\rm (H5)}{=}\mathbb{I}\Big(N, \mathbb{H}^{-1}(Q)\wedge \mathbb{H}^{-1}\mathbb{H}^{-1}(Q)\Big)\stackrel{\rm (HM)}{=}\mathbb{I}\Big(Q,\mathbb{H}(N)\Big),$$ i.e., (HSM) holds.

(ii) (HSM) $\Longrightarrow$ (HM)+(H5).

Analogizing the verification  of (i) in Theorem \ref{osutheoremr}, we have that (HSM) $\Longrightarrow$ (H0). Then $\forall b\in X$,  $E\in P(\mathbb{U})$  $$\mathbb{H}^{-1}(E)(b)\stackrel{\rm Proposition\   \ref{lemmafn}(3)}{=}\mathbb{I}\Big(E,\mathbb{H}(\mathbb{U}_{\{b\}})\Big)\stackrel{\rm (HSM)}{=}\mathbb{I}\Big(\mathbb{U}_{\{b\}}, \mathbb{H}\mathbb{H}(E)\wedge \mathbb{H}(E)\Big)\leq  \mathbb{H}(E)(b),$$ similarly $ \mathbb{H}^{-1}(E)(b)\geq  \mathbb{H}(E)(b)$, so (H5) gains. Using (H5) to (HSM) one obtain (HM).

\end{proof}

Next, we consider the combinations of three conditions.

\begin{theorem} \label{osutheoremrst}  $\mathbb{H}$ is a reflexive, transitive and symmetric, i.e., equivalent  {\it LVHRAO} iff it fulfills (HRTS).

(HRTS) $\forall Q, N\in P(\mathbb{U})$,
$\mathbb{I}\Big(Q,\mathbb{H}(N)\Big)=\mathbb{I}\Big(N,  Q\vee \mathbb{H}(Q)\vee \mathbb{H}\mathbb{H}(Q)\Big).$
\end{theorem}

\begin{proof}  From Theorem \ref{osutheoremrt} and Proposition \ref{oprops} we only need examine (HRT)+(H5) $\Longleftrightarrow$ (HRTS).

(i) (HRT)+(H5) $\Longrightarrow$ (HRTS). Indeed $$\mathbb{I}\Big(N,  Q\vee \mathbb{H}(Q)\vee \mathbb{H}\mathbb{H}(Q)\Big)\stackrel{\rm (H5)}{=}\mathbb{I}\Big(N, Q\vee \mathbb{H}^{-1}(Q)\vee \mathbb{H}^{-1}\mathbb{H}^{-1}(Q)\Big)\stackrel{\rm (HRT)}{=}\mathbb{I}\Big(Q,\mathbb{H}(N)\Big),$$ i.e., (HRTS) holds.

(ii) (HRTS) $\Longrightarrow$ (HRT)+(H5).

Analogizing the verification  of (i) in Theorem \ref{osutheoremr}, we have that (HRTS) $\Longrightarrow$ (H0). Then $\forall b\in X$,  $K\in P(\mathbb{U})$ $$\mathbb{H}^{-1}(K)(b)\stackrel{\rm Proposition\   \ref{lemmafn}(3)}{=}\mathbb{I}\Big(K,\mathbb{H}(\mathbb{U}_{\{b\}})\Big)\stackrel{\rm (HRTS)}{=}\mathbb{I}\Big(\mathbb{U}_{\{b\}}, K\vee \mathbb{H}\mathbb{H}(K)\vee \mathbb{H}(K)\Big)\geq  \mathbb{H}(K)(b),$$ similarly $ \mathbb{H}^{-1}(K)(b)\leq  \mathbb{H}(K)(b)$, so (H5) gains. Using (H5) to (HRTS) one get (HRT).
\end{proof}

\begin{proposition} \label{osutheoremrte}

(1) $\mathbb{H}$ is a reflexive, transitive and Euclidean {\it LVHRAO} iff it fulfills (HRTE).

(HRTE) $\forall Q, N\in P(\mathbb{U})$, $\mathbb{I}\Big(Q,\mathbb{H}(N)\Big)=\mathbb{I}\Big(N, Q\vee \mathbb{H}^{-1}\mathbb{H}^{-1}(Q)\vee \mathbb{H}\mathbb{H}^{-1}(Q)\vee \mathbb{H}^{-1}(Q)\Big).$

(2) $\mathbb{H}$ is a reflexive, symmetric and Euclidean  {\it LVHRAO} iff it fulfills (HRTS).

(3)  $\mathbb{H}$ is a transitive, symmetric and Euclidean  {\it LVHRAO} iff it fulfills (HTS).

\end{proposition}

\begin{proof} (1) By Theorem \ref{osutheoremrt} and Proposition \ref{opropeu}, we only need examine (HRT)+(H6) $\Longleftrightarrow$ (HRTE). The detailed proof is similar to Theorem \ref{osutheoremrst}.

(2) It follows by Theorem \ref{osutheoremrst} and that if $\mathbb{R}$ is symmetric, then $\mathbb{R}$ is transitive iff $\mathbb{R}$ is Euclidean.

(3) It follows by Proposition \ref{osutheoremts} (2) and that if $\mathbb{R}$ is symmetric, then $\mathbb{R}$ is transitive iff $\mathbb{R}$ is Euclidean.

\end{proof}

At last, we give the least (resp., largest) equivalent {\it LVHRAO} on $\mathbb{U}$ and the corresponding {\it Lvr}.

\begin{example} (1) Let $\mathbb{H}_0: P(\mathbb{U})\longrightarrow P(\mathbb{U})$ be defined by $\forall Q\in P(\mathbb{U}), \mathbb{H}_0(Q)=Q$. Then $\forall Q, N\in P(\mathbb{U})$,
$$\mathbb{I}\Big(Q,\mathbb{H}_0(N)\Big)=\mathbb{I}\Big(Q,N\Big)=\mathbb{I}\Big(N,  Q\Big)=\mathbb{I}\Big(N,  Q\vee \mathbb{H}_0(Q)\vee \mathbb{H}_0\mathbb{H}_0(Q)\Big),$$ i.e., (HRTS) holds.  According to Theorem \ref{osutheoremrst}, there is an uniquely  equivalent {\it Lvr} $\mathbb{R}$ on $\mathbb{U}$  s.t. $\overline{\mathbb{R}_0}=\mathbb{H}_0$. Hence, $\forall a,b \in X$, $$\mathbb{R}_0(a,b)=\overline{\mathbb{R}_0}(\mathbb{U}_{\{a\}})(b)=\mathbb{H}_0(\mathbb{U}_{\{a\}})(b)=\mathbb{U}_{\{a\}}(b),$$ which means $\mathbb{R}_0$ is the least equivalent {\it Lvr} on $\mathbb{U}$, and $\mathbb{H}_0$ is the least equivalent {\it LVHRAO} on $\mathbb{U}$.

(2)  Let $\mathbb{H}_1: P(\mathbb{U})\longrightarrow P(\mathbb{U})$ be defined by $$\forall Q\in P(\mathbb{U}), \forall a\in X,  \mathbb{H}_1(Q)(a)=\bigvee_{b\in X}\Big(Q(b)\circledast(\mathbb{U}(b)\rightarrow \mathbb{U}(a))\Big).$$ Then $$\mathbb{H}_1(Q)(a)=\bigvee_{b\in X}\Big(Q(b)\circledast(\mathbb{U}(b)\rightarrow \mathbb{U}(a))\Big)\geq Q(a)\circledast(\mathbb{U}(a)\rightarrow \mathbb{U}(a))=Q(a),$$ and $$\mathbb{H}_1\mathbb{H}_1(Q)(a)=\bigvee_{b,c\in X}\Big(Q(c)\circledast (\mathbb{U}(c)\rightarrow \mathbb{U}(b))\circledast(\mathbb{U}(b)\rightarrow \mathbb{U}(a))\Big)\leq\bigvee_{c\in X}\Big(Q(c)\circledast(\mathbb{U}(c)\rightarrow \mathbb{U}(a))\Big)=\mathbb{H}_1(Q)(a).$$It follows that $\forall Q, N\in P(\mathbb{U})$, \begin{eqnarray*}\mathbb{I}\Big(N,  Q\vee \mathbb{H}_1(Q)\vee \mathbb{H}_1\mathbb{H}_1(Q)\Big)&=&\mathbb{I}\Big(N,  \mathbb{H}_1(Q)\Big)\\
&=&\bigvee_{a,b\in X}\Big(Q(a)\circledast (\mathbb{U}(a)\rightarrow \mathbb{U}(b))\circledast(\mathbb{U}(b)\rightarrow N(b))\Big), \ {\rm by \ Lemma \ref{jb}(1)}\\
&=&\bigvee_{a,b\in X}\Big(\mathbb{U}(a)\circledast(\mathbb{U}(a)\rightarrow Q(a))\circledast (\mathbb{U}(a)\rightarrow \mathbb{U}(b))\circledast(\mathbb{U}(b)\rightarrow N(b))\Big)\\
&=&\bigvee_{a,b\in X}\Big((\mathbb{U}(a)\rightarrow Q(a))\circledast (\mathbb{U}(a)\wedge \mathbb{U}(b))\circledast(\mathbb{U}(b)\rightarrow N(b))\Big), \ {\rm by \ Lemma \ref{jb}(2)}\\
&=&\bigvee_{a,b\in X}\Big((\mathbb{U}(a)\rightarrow Q(a))\circledast N(b)\circledast(\mathbb{U}(b)\rightarrow(\mathbb{U}(a)\wedge \mathbb{U}(b)))\Big)\\
&=&\bigvee_{a,b\in X}\Big(N(b)\circledast(\mathbb{U}(b)\rightarrow\mathbb{U}(a))\circledast(\mathbb{U}(a)\rightarrow Q(a))\Big)\\
&=&\mathbb{I}\Big(Q,\mathbb{H}_1(N)\Big),
\end{eqnarray*} i.e., (HRTS) holds.  According to Theorem \ref{osutheoremrst}, there is an uniquely  equivalent {\it Lvr} $\mathbb{R}_1$ on $\mathbb{U}$  s.t. $\overline{\mathbb{R}_1}=\mathbb{H}_1$. Hence, $\forall a,b \in X$, $$\mathbb{R}_1(a,b)=\overline{\mathbb{R}_1}(\mathbb{U}_{\{a\}})(b)=\mathbb{H}_1(\mathbb{U}_{\{a\}})(b)=\bigvee_{c\in X}\Big(\mathbb{U}_{\{a\}}(c)\circledast(\mathbb{U}(c)\rightarrow \mathbb{U}(b))\Big)=\mathbb{U}(a)\wedge \mathbb{U}(b),$$ which means $\mathbb{R}_1$ is the largest equivalent {\it Lvr} on $\mathbb{U}$, and $\mathbb{H}_1$ is the largest equivalent {\it LVHRAO} on $\mathbb{U}$.
\end{example}

\section{A novel axiomatic approach to {\it LVLRAO}s}

This section  provides the notion of outer product (and  the subsethood degree) of $L$-subsets in $\mathbb{U}$, and applies it in novel  axiomatization  on {\it LVLRAO} and  some special {\it LVLRAO}s associated with reflexive, symmetric, transitive, Euclidean and mediate {\it Lvr}s and their compositions.

In this section, following reference \cite{LF19,BP21}, we assume that $(L,\circledast)$ is a complete MV-algebra, that is, a GL-quantale $(L,\circledast)$ with $\forall \gamma\in L$, $(\gamma\rightarrow 0)\rightarrow 0=\gamma$.

\begin{lemma} \cite{BP21} \label{jbmv} $(L,\circledast)$ be a complete MV-algebra and $\alpha_{i}(i\in I), \beta,\gamma,\in L$, $W, V\in P(\mathbb{U})$ and $b\in X$.

{\rm (1)} $\beta\rightarrow \gamma=(\gamma\rightarrow 0)\rightarrow (\beta\rightarrow 0)$.

{\rm (2)}  $\beta\circledast \bigwedge\limits_{i}\alpha_{i}=\bigwedge\limits_{i}(\beta\circledast \alpha_{i})$.

{\rm (3)} $W\rightarrow V =\neg V \rightarrow \neg W $, $\neg \neg W =W$.

{\rm (4)} $\neg\mathbb{U}_{X-\{b\}}=\mathbb{U}_{\{b\}}$.

{\rm (5)} $\mathbb{U}(b)\rightarrow \neg W(b)=W(b)\rightarrow 0$.
\end{lemma}

\subsection{The outer product of $L$-subsets in $\mathbb{U}$ and the lower inverse mapping of $\mathbb{L}: P(\mathbb{U})\longrightarrow P(\mathbb{U})$}


\begin{definition} (1) Define mapping $\mathbb{S}: P(\mathbb{U})\times P(\mathbb{U})\longrightarrow L$  by $\forall M, Q\in P(\mathbb{U})$:
$$\mathbb{S}(M,Q)=\bigwedge\limits_{b\in X}\Big(\mathbb{U}(b)\circledast \big(M(b)\rightarrow  Q(b)\big)\Big).$$  Then  $ \mathbb{S}(M,Q)$ is called the  subsethood degree of $M,Q$.

(2) Define mapping $\mathbb{O}: P(\mathbb{U})\times P(\mathbb{U})\longrightarrow L$  by $\forall M, Q\in P(\mathbb{U})$:
$$\mathbb{O}(M,Q)=\mathbb{S}(\neg M,Q)=\bigwedge\limits_{b\in X}\Big(\mathbb{U}(b)\circledast \big(\neg M(b)\rightarrow  Q(b)\big)\Big).$$  Then  $ \mathbb{O}(M,Q)$ is called the  outer product of $M,Q$.
\end{definition}

\begin{definition}  The mapping  $\mathbb{L}^{\sim 1}: P(\mathbb{U})\longrightarrow   P(\mathbb{U})$ associated with $\mathbb{L}: P(\mathbb{U})\longrightarrow   P(\mathbb{U})$ is called the lower inverse mapping of $\mathbb{L}$, where $\forall Q\in  P(\mathbb{U}), b\in X$: $ \mathbb{L}^{\sim 1}(Q)(b)=\mathbb{O}\Big(\mathbb{L}(\mathbb{U}_{X-\{b\}}), Q\Big).$
\end{definition}

\begin{remark} Let $\mathbb{U}=1$ and $M,Q\in P(\mathbb{U})$.

(1) Note that $$\mathbb{S}(M,Q)=\bigwedge\limits_{b\in X}\Big(\mathbb{U}(b)\circledast \big(M(b)\rightarrow  Q(b)\big)\Big)=\bigwedge\limits_{b\in X}\Big(M(b)\rightarrow Q(b)\Big).$$Hence $\mathbb{S}(M,Q)$ degenerates into the subsethood degree in  \cite{RB12}.

(2) Note that $$\mathbb{O}(M,Q)=\bigwedge\limits_{b\in X}\Big(\mathbb{U}(b)\circledast \big(\neg M(b)\rightarrow  Q(b)\big)\Big)=\bigwedge\limits_{b\in X}\Big((M(b)\rightarrow 0)\rightarrow Q(b)\Big).$$Hence $\mathbb{O}(M,Q)$ degenerates into the outer product in  \cite{Bao18}.

(3) Note that
$$\mathbb{L}^{\sim 1}(M)(d)=\mathbb{O}\Big( \mathbb{L}(\mathbb{U}_{X-\{d\}}), M\Big)=\bigwedge\limits_{a\in X}\Big((\mathbb{L}(\chi_{X-\{d\}})(a)\rightarrow 0)\rightarrow M(a)\Big).$$
Hence $\mathbb{L}^{\sim 1}$  degenerates into the lower inverse mapping in \cite{Bao18}.
\end{remark}

\begin{example} Let $X=\{a,b,c,d\}$, $\mathbb{U}=\frac{0.2}{a}+\frac{0.7}{b}+\frac{0.3}{c}+\frac{0.8}{d}$ and $\mathbb{L}: P(\mathbb{U})\longrightarrow P(\mathbb{U})$ be defined by $\forall Q\in P(\mathbb{U}), \mathbb{L}(Q)=Q$.  Consider $L=[0,1]$ and $\circledast$ being defined through $\forall \alpha, \beta\in [0,1]$, $\alpha\circledast \beta=\max\{\alpha+\beta-1, 0\}$,  then $([0,1],\circledast)$ forms a complete MV-algebra.

Take $M=\frac{0.2}{a}+\frac{0.5}{b}+\frac{0.3}{c}+\frac{0.6}{d}$ and $Q=\frac{0.2}{a}+\frac{0.5}{b}+\frac{0.3}{c}+\frac{0.5}{d}$, then $M, Q\in P(\mathbb{U})$. So, $$\mathbb{O}(M,Q)=0.2\wedge 0.7\wedge 0.3 \wedge 0.8=0.2,  \mathbb{L}^{\sim 1}(Q)=\frac{0.2}{a}+\frac{0.2}{b}+\frac{0.2}{c}+\frac{0.2}{d}.$$
\end{example}

After careful analysis, we find that for any $L$-universe $\mathbb{U}$, we can't prove that $\mathbb{O}$ and $\mathbb{L}^{\sim 1}$ have the expected properties. But when $\mathbb{U}$ is a constant $L$-universe (i.e., $\mathbb{U}(b)=\mathbb{U}(c), \forall b,c \in X$), the goal can be achieved. Therefore, in the following,  we always assume that $\mathbb{U}$ is a constant $L$-universe.

\begin{proposition} \label{lemmaOIL} Let $Q, P, P_i(i\in I)\in P(\mathbb{U})$.

{\rm (1)} $\forall D\in P(\mathbb{U}), \mathbb{O}\Big(D, Q\Big)\leq \mathbb{O}\Big(D, P\Big)$ implies $Q\leq P$, furthermore, $\mathbb{O}\Big(D, Q\Big)= \mathbb{O}\Big(D, P\Big)$ implies $Q= P$.

{\rm (2)} $\mathbb{O}\Big(P, Q\Big)=\mathbb{O}\Big(Q, P\Big)$.

{\rm (3)} $\forall \beta\in  L, b\in X$,
$\mathbb{O}\Big(K,  \mathbb{U}\wedge (\beta\rightarrow P)\Big)=\mathbb{U}(b)\wedge\Big(\beta\rightarrow \mathbb{O}\Big(K, P\Big)\Big).$

{\rm (4)} $\mathbb{O}\Big(Q,  \bigwedge\limits_{i} P_i\Big)=\bigwedge\limits_{i}\mathbb{O}\Big(Q, P_i\Big)$.

{\rm (5)}  $\forall b, h \in X$, $\mathbb{L}^{\sim 1} (\mathbb{U}_{X-\{h\}})(b)=\mathbb{L} (\mathbb{U}_{X-\{b\}})(h).$
\end{proposition}

\begin{proof}  (1) $\forall b\in X$,
\begin{eqnarray*}\mathbb{O}\Big(\mathbb{U}_{X-\{b\}}, Q\Big)&=&\bigwedge_{z\in X}\Big(\mathbb{U}(z)\circledast \big(\neg\mathbb{U}_{X-\{b\}}(z)\rightarrow Q(z)\big)\Big)\\
&=& \Big(\bigwedge_{z\neq b}\Big(\mathbb{U}(z)\circledast \big(\neg\mathbb{U}_{X-\{b\}}(z)\rightarrow Q(z)\big)\Big)\Big)\wedge \Big(
\mathbb{U}(b)\circledast \big(\mathbb{U}(b)\rightarrow Q(b)\big)\Big)\\
&=&\Big(\bigwedge_{z\neq b}\mathbb{U}(z)\Big)\wedge Q(b)\ {\rm by }\  \mathbb{U}(z)=\mathbb{U}(b)\\
&=& Q(b).\end{eqnarray*} Hence $\mathbb{O}\Big(\mathbb{U}_{X-\{b\}}, Q\Big)\leq \mathbb{O}\Big(\mathbb{U}_{X-\{b\}}, P)$ implies $Q\leq P$. The further part follows similar.

(2) It is concluded by Lemma \ref{jbmv}(3).

(3) For any $d\in X, \beta\in L$,
\begin{eqnarray*}\mathbb{O}\Big(K,  \mathbb{U}\wedge (\beta\rightarrow P)\Big)&=&\bigwedge\limits_{z\in X}\Big(\mathbb{U}(z)\circledast \big(\neg K(z)\rightarrow  \mathbb{U}(z)\wedge (\beta\rightarrow P(z))\big)\Big)\\
&=&\bigwedge\limits_{z\in X}\Big(\mathbb{U}(z)\wedge\big( (\mathbb{U}(z)\rightarrow \neg K(z))\rightarrow  (\beta\rightarrow P(z))\big)\Big)\\
&=&\bigwedge\limits_{z\in X}\Big(\mathbb{U}(z)\wedge\big( \beta\rightarrow  ( (\mathbb{U}(z)\rightarrow \neg K(z))\rightarrow P(z))\big)\Big)\\
&=&\bigwedge\limits_{z\in X}\Big(\mathbb{U}(z)\wedge\big( \beta\rightarrow  \mathbb{U}(z)\wedge( (\mathbb{U}(z)\rightarrow \neg K(z))\rightarrow P(z))\big)\Big)\ {\rm by }\  \mathbb{U}(z)=\mathbb{U}(d)\\
&=&\mathbb{U}(d) \wedge \Big(\beta\rightarrow\bigwedge\limits_{z\in X}[\mathbb{U}(z)\circledast (\neg K(z)\rightarrow  P(z))]\Big)\\
&=&\mathbb{U}(d) \wedge \Big(\beta\rightarrow\mathbb{O}\Big(K,P\Big)\Big).\end{eqnarray*}

(4) It holds by $\mathbb{O}\Big(Q,  \bigwedge\limits_{i} P_i\Big)= \bigwedge\limits_{a\in X}\bigwedge\limits_{i}\Big(\mathbb{U}(a)\circledast \big(\neg Q(a)\rightarrow P_i(a)\big)\Big)=\bigwedge\limits_{i}\mathbb{O}\Big(Q, P_i\Big).$

 (5) It holds by
\begin{eqnarray*} \mathbb{L}^{\sim 1}\Big(\mathbb{U}_{X-\{h\}}\Big)(b)&=&\bigwedge_{g\in X} \Big(\mathbb{U}(g)\circledast \big(\neg \mathbb{L}(\mathbb{U}_{X-\{b\}})(g)\rightarrow\mathbb{U}_{X-\{h\}}(g)\big)\Big)\\
&=& \mathbb{U}(g)\wedge \Big( \mathbb{U}(h)\circledast \big(\neg \mathbb{L}(\mathbb{U}_{X-\{b\}})(h)\rightarrow 0\big)\Big)\ {\rm by }\  \mathbb{U}(g)=\mathbb{U}(h)\\
&=& \mathbb{U}(h)\wedge \Big( \mathbb{U}(h)\circledast \big( \mathbb{U}(h)\rightarrow \mathbb{L}(\mathbb{U}_{X-\{b\}})(h)\big)\Big)\\
&=&\mathbb{L}(\mathbb{U}_{X-\{b\}})(h).\qedhere\end{eqnarray*} \end{proof}

\subsection{A novel single axiom characterizations about {\it LVLRAO}s}

This subsection  provides the single axiom characterizations about {\it LVLRAO} and  some special {\it LVLRAO}s associated with reflexive, symmetric, transitive, Euclidean $\&$ mediate {\it Lvr}s, respectively.

\begin{lemma} \label{lemmafenl}  For  $Q\in  P(\mathbb{U})$, $Q=\bigwedge\limits_{h\in X}\Big(\mathbb{U}\wedge \big ((Q(h)\rightarrow 0)\rightarrow \mathbb{U}_{X-\{h\}}\big)\Big)$.
\end{lemma}

\begin{proof} That is concluded from  $\forall b\in X$, \begin{eqnarray*}& &\bigwedge\limits_{h\in X}\Big(\mathbb{U}\wedge \big ((Q(h)\rightarrow 0)\rightarrow \mathbb{U}_{X-\{h\}}\big)\Big)(b)\\
&=&\bigwedge\limits_{h\neq b}\Big(\mathbb{U}(b)\wedge \big ((Q(h)\rightarrow 0)\rightarrow \mathbb{U}_{X-\{h\}}(b)\big)\Big)\wedge \Big(\mathbb{U}(b)\wedge \big((Q(b)\rightarrow 0)\rightarrow 0\big)\Big)\\
&=&\Big(\mathbb{U}(b)\wedge \big ((Q(h)\rightarrow 0)\rightarrow \mathbb{U}(b)\big)\Big)\wedge Q(b)\\
&=&Q(b).\qedhere\end{eqnarray*}\end{proof}

Unless otherwise stated, we presume  $\mathbb{L}: P(\mathbb{U})\longrightarrow   P(\mathbb{U})$ to be a mapping. We introduce the following the following notations.

(L1)  $\mathbb{L}(\mathbb{U}\wedge (\alpha\rightarrow W))=\mathbb{U}\wedge (\alpha\rightarrow \mathbb{L} (W))$ for any $\alpha\in L, W \in P(\mathbb{U})$.

(L2)  $\mathbb{L}   (\bigwedge\limits_{i}W_i)=\bigwedge\limits_{i} \mathbb{L}   (W_i)$ for any $W_i(i\in I)\subseteq P(\mathbb{U})$.

\begin{remark} \label{remgerl} From Proposition \ref{theoremL1} (3) we know that $\mathbb{L}$ is an {\it LVLRAO} iff it fulfills (L1) and (L2).\end{remark}

\begin{theorem} \label{sutheorembgl}  $\mathbb{L}$ is an LVLRAO iff it fulfills (L).

(L) $ \forall  M, Q \in P(\mathbb{U})$,
$ \mathbb{O}\Big(Q, \mathbb{L} (M)\Big)=\mathbb{O}\Big(M,  \mathbb{L}^{\sim 1}(Q)\Big).$
\end{theorem}

\begin{proof} From Remark \ref{remgerl} we only to examine   (L1)+(L2) $\Longleftrightarrow$ (L).

(1) (L1), (L2) $\Longrightarrow$ (L). $\forall M, Q\in P(\mathbb{U})$,  \begin{eqnarray*}& &\mathbb{O}\Big(M,  \mathbb{L}^{\sim 1}(Q)\Big)\\
&=&\bigwedge_{b\in X} \Big(\mathbb{U}(b)\circledast\big(\neg M(b)\rightarrow \mathbb{L}^{\sim 1}(Q)(b)\big)\Big)   \ {\rm by \ Lemma\   \ref{jb}}(3)\\
&=&\bigwedge_{b\in X} \Big(\mathbb{U}(b)\wedge \Big[(\mathbb{U}(b)\rightarrow \neg M(b))\rightarrow \mathbb{L}^{\sim 1}(Q)(b)\Big]\Big)   \ {\rm by \ Lemma\   \ref{jbmv}}(5)\\
&=&\bigwedge_{b\in X} \Big(\mathbb{U}(b)\wedge \Big[(M(b)\rightarrow 0)\rightarrow \bigwedge_{z\in X} \Big\{\mathbb{U}(z)\wedge [(Q(z)\rightarrow 0)\rightarrow \mathbb{L}(\mathbb{U}_{X-\{b\}})(z)]\Big\}\Big]\Big), \mathbb{U}(b)=\mathbb{U}(z)\\
&=&\bigwedge_{b,z\in X} \Big(\mathbb{U}(z)\wedge \Big[(M(b)\rightarrow 0)\rightarrow \Big\{\mathbb{U}(z)\wedge [(Q(z)\rightarrow 0)\rightarrow \mathbb{L}(\mathbb{U}_{X-\{b\}})(z)]\Big\}\Big]\Big)\\
&=&\bigwedge_{b,z\in X} \Big(\mathbb{U}(z)\wedge \Big[(M(b)\rightarrow 0)\rightarrow [(Q(z)\rightarrow 0)\rightarrow \mathbb{L}(\mathbb{U}_{X-\{b\}})(z)]\Big]\Big)\\
&=&\bigwedge_{z\in X} \Big(\mathbb{U}(z)\wedge \Big[(Q(z)\rightarrow 0)\rightarrow \bigwedge_{b\in X} [(M(b)\rightarrow 0)\rightarrow \mathbb{L}(\mathbb{U}_{X-\{b\}})(z)]\Big]\Big)\\
&=&\bigwedge_{z\in X} \Big(\mathbb{U}(z)\wedge \Big[(Q(z)\rightarrow 0)\rightarrow \bigwedge_{b\in X} \Big\{\mathbb{U}(z)\wedge [(M(b)\rightarrow 0)\rightarrow \mathbb{L}(\mathbb{U}_{X-\{b\}})(z)]\Big\}\Big]\Big)  \ {\rm by \ (H1),(H2)}\\
&=&\bigwedge_{z\in X} \Big(\mathbb{U}(z)\wedge \Big[(Q(z)\rightarrow 0)\rightarrow \mathbb{L}\Big(\bigwedge_{b\in X}\mathbb{U}\wedge [(M(b)\rightarrow 0)\rightarrow \mathbb{U}_{X-\{b\}}]\Big)(z)\Big]\Big)   \ {\rm by \ Lemma\   \ref{lemmafenl}}\\
&=&\bigwedge_{z\in X} \Big(\mathbb{U}(z)\wedge \Big[(Q(z)\rightarrow 0)\rightarrow \mathbb{L}(M)(z)\Big]\Big)\\
&=&\bigwedge_{z\in X} \Big(\mathbb{U}(z)\wedge \Big[(\mathbb{U}(z)\rightarrow \neg Q(z))\rightarrow \mathbb{L}(M)(z)\Big]\Big)\\
&=&\bigwedge_{z\in X} \Big(\mathbb{U}(z)\circledast \Big[\neg Q(z)\rightarrow \mathbb{L}(M)(z)\Big]\Big)=\mathbb{O}\Big(Q, \mathbb{L} (M)\Big).\end{eqnarray*}

(2) (L)$\Longrightarrow$ (L1), (L2).

(i) Let $\beta\in  L, V\in  P(\mathbb{U})$, $d\in X$. Then

\begin{eqnarray*}\mathbb{L} \big(\mathbb{U}\wedge (\beta\rightarrow V)\big)(d)&=&\mathbb{O}\Big(\mathbb{U}_{X-\{d\}},  \mathbb{L} \big(\mathbb{U}\wedge (\beta\rightarrow V)\big)\Big)\\
&\stackrel{\rm (L)}{=}& \mathbb{O}\Big(\mathbb{U}\wedge(\beta\rightarrow V),   \mathbb{L}^{\sim 1}(\mathbb{U}_{X-\{d\}})\Big)\ {\rm by \ Proposition\   \ref{lemmaOIL} (2),(3)}
\\
&=&\mathbb{U}(d)\wedge\Big(\beta\rightarrow \mathbb{O}\Big(V, \mathbb{L}^{\sim 1}(\mathbb{U}_{X-\{d\}})\Big)\Big)\\
&\stackrel{\rm (L)}{=}&\mathbb{U}(d)\wedge\Big(\beta\rightarrow \mathbb{O}\Big(\mathbb{U}_{X-\{d\}}, \mathbb{L} (V)\Big)= \mathbb{U}(d)\wedge\Big(\beta\rightarrow \mathbb{L} (V)(d)\Big).\end{eqnarray*}Hence $ \mathbb{L} (\mathbb{U}\wedge (\beta\rightarrow V))= \mathbb{U}\wedge (\beta\rightarrow \mathbb{L} (V))$, i.e., (L1) holds.

(ii) $\forall  Q, V_i(i\in I)\in  P(\mathbb{U})$, by (L) and  Proposition \ref{lemmaOIL} (2),(4)
\begin{eqnarray*}\mathbb{O}\Big(Q, \mathbb{L} (\bigwedge_{i} V_i)\Big)&=& \mathbb{O}\Big(  \bigwedge_{i} V_i,   \mathbb{L}^{\sim 1}(Q)\Big)=\bigwedge_{i}\mathbb{O}\Big(\mathbb{L}^{\sim 1}(Q), V_i\Big)=\mathbb{O}\Big(Q, \bigwedge_{i} \mathbb{L} (V_i)\Big).\end{eqnarray*}This shows   $ \mathbb{L} (\bigwedge\limits_{i} V_i)= \bigwedge\limits_{i}  \mathbb{L} ( V_i)$, i.e., (L2) holds.
\end{proof}

To facilitate the expression, we give the following concepts.

\begin{definition} A mapping $\mathbb{L}: P(\mathbb{U})\longrightarrow P(\mathbb{U})$ is named a reflexive {\it LVLRAO} if there exists a reflexive {\it Lvr} $\mathbb{R}$ on $\mathbb{U}$ s.t. $\mathbb{L}=\underline{\mathbb{R}}$. Similarly, one can define many other special  {\it LVLRAO}.

\end{definition}

\begin{remark} Next, we will describe some special {\it LVLRAO}s by strengthening axiom (L). For readers to quickly understand the operators described by these axioms, we will abbreviate ``reflexive, transitive, symmetric, Euclidean and mediate" as ``R, T, S, E and M" respectively. Then, readers can understand that axiom (LR)  describes reflexive {\it LVLRAO}, (LRT) describes reflexive and transitive {\it LVLRAO}, (LRTS) describes reflexive, transitive and symmetric {\it LVLRAO}.
\end{remark}

\begin{proposition} \label{oproprl}  $\mathbb{L}$ is a reflexive LVLRAO iff it fulfills (L1), (L2) and (L3).

(L3) $\forall W\in  P(\mathbb{U})$, $\mathbb{L}^{\sim 1}(W)\leq W$.
\end{proposition}

\begin{proof}
  $\Longrightarrow$. Let $\mathbb{L} =\underline{\mathbb{R}}$ for a reflexive {\it Lvr} $\mathbb{R}$ on $\mathbb{U}$. We only examine (L3). $\forall W\in  P(\mathbb{U})$, $b\in X$, \begin{eqnarray*}\mathbb{L}^{\sim 1}(W)(b)&=&\bigwedge\limits_{z\in X}\Big(\mathbb{U}(z)\circledast (\neg\underline{\mathbb{R}}(\mathbb{U}_{X-\{b\}})(z)\rightarrow W(z)) \Big)\\
 &=&\bigwedge\limits_{z\in X}\Big(\mathbb{U}(z)\circledast (\mathbb{R}(b,z)\rightarrow W(z)) \Big)\\
 &\leq&\mathbb{U}(b)\circledast (\mathbb{R}(b,b)\rightarrow W(b)) \ {\rm by \ reflexivity}\\
 &\leq& \mathbb{U}(b)\circledast (\mathbb{U}(b)\rightarrow W(b))= W(b).\end{eqnarray*}

$\Longleftarrow$. By Remark \ref{remgerl}, $\mathbb{L} =\underline{\mathbb{R}}$ for an {\it Lvr} $\mathbb{R}$. Then $\forall b\in X$,  \begin{eqnarray*}\mathbb{U}(b)&=&\neg\mathbb{U}_{X-\{b\}}(b)\\
&\stackrel{(L3)}{\leq}& \neg \mathbb{L}^{\sim 1}(\mathbb{U}_{X-\{b\}})(b)   \ {\rm by \ Proposition\   \ref{lemmaOIL}}(5)\\
&=& \neg\underline{\mathbb{R}}(\mathbb{U}_{X-\{b\}})(b)=\mathbb{R}(b,b).\end{eqnarray*} So, $\mathbb{R}$ is reflexive.
\end{proof}

\begin{theorem} \label{osutheoremrl} $\mathbb{L}$ is a reflexive {\it LVLRAO} iff it fulfills (LR).

(LR) $\forall M, Q\in P(\mathbb{U})$, $\mathbb{O}\Big(Q,\mathbb{L}(M)\Big)=\mathbb{O}\Big(M, Q\wedge \mathbb{L}^{\sim 1}(Q)\Big).$
\end{theorem}

\begin{proof}  By Theorem \ref{sutheorembgl} and Proposition \ref{oproprl}, we only need examine (L)+(L3) $\Longleftrightarrow$ (LR).

(1) (L)+(L3) $\Longrightarrow$ (LR). Indeed  $$\mathbb{O}\Big(M, Q\wedge \mathbb{L}^{\sim 1}(Q)\Big)\stackrel{\rm (L3)}{=}\mathbb{O}\Big(M, \mathbb{L}^{\sim 1}(Q)\Big)\stackrel{\rm (L)}{=}\mathbb{O}\Big(Q,\mathbb{L}(M)\Big),$$ i.e., (LR) holds.

(2) (LR) $\Longrightarrow$ (L)+(L3).$\forall o\in X$,  put $M=\mathbb{U}_{X-\{o\}}$ in (LR) then
$$\mathbb{L}^{\sim 1}(Q)(o)\stackrel{\rm Proposition\   \ref{lemmaOIL}(2)}{=}\mathbb{O}\Big(Q,\mathbb{L}(\mathbb{U}_{X-\{o\}})\Big)\stackrel{\rm (LR)}{=}\mathbb{O}\Big(\mathbb{U}_{X-\{o\}}, Q\wedge \mathbb{L}^{\sim 1}(Q)\Big)\leq  Q(o),$$ so (L3) gains. Using (L3) to (LR) one get (L).
\end{proof}

\begin{proposition} \label{oproptl}  $\mathbb{L}$ is a transitive LVLRAO iff it fulfills (L1), (L2) and (L4).

 (L4) $\forall W \in P(\mathbb{U})$, $\mathbb{L}^{\sim 1}\mathbb{L}^{\sim 1}(W)\geq \mathbb{L}^{\sim 1} (W)$.
\end{proposition}

\begin{proof} $\Longrightarrow$.  Let $\mathbb{L} =\underline{\mathbb{R}}$ for a transitive {\it Lvr} $\mathbb{R}$ on $\mathbb{U}$. We only examine (L4). Indeed, $\forall W\in  P(\mathbb{U})$, $d\in X$,
 \begin{eqnarray*}\mathbb{L}^{\sim 1}\mathbb{L}^{\sim 1}(W)(d)&=&\bigwedge_{b\in X}\Big(\mathbb{U}(b)\circledast \Big[\neg\underline{\mathbb{R}}(\mathbb{U}_{X-\{d\}})(b)\rightarrow \underline{\mathbb{R}}^{\sim 1} (W)(b)\Big]\Big)\\
&=&\bigwedge_{b\in X}\Big(\mathbb{U}(b)\circledast \Big[\mathbb{R}(d,b)\rightarrow \bigwedge_{g\in X}\Big\{\mathbb{U}(g)\circledast [\neg\underline{\mathbb{R}}(\mathbb{U}_{X-\{b\}})(g)\rightarrow W(g)]\Big\}\Big]\Big)   \ {\rm by}\  \mathbb{U}(b)=\mathbb{U}(g)\\
&=&\bigwedge_{b\in X}\Big(\mathbb{U}(g)\circledast \Big[\mathbb{R}(d,b)\rightarrow \bigwedge_{g\in X}\Big\{\mathbb{U}(b)\circledast [\mathbb{R}(b,g)\rightarrow W(g)]\Big\}\Big]\Big)\\
&=&\bigwedge_{b,g\in X}\Big(\mathbb{U}(g)\circledast \Big[\mathbb{R}(d,b)\rightarrow \Big\{\mathbb{U}(b)\wedge [(\mathbb{U}(b)\rightarrow \mathbb{R}(b,g))\rightarrow W(g)]\Big\}\Big]\Big)\\
&=&\bigwedge_{b,g\in X}\Big(\mathbb{U}(g)\circledast \Big[\mathbb{R}(d,b)\rightarrow  [(\mathbb{U}(b)\rightarrow \mathbb{R}(b,g))\rightarrow W(g)]\Big]\Big)\\
&=&\bigwedge_{b,g\in X}\Big(\mathbb{U}(g)\circledast \Big[[\mathbb{R}(d,b)\circledast(\mathbb{U}(b)\rightarrow \mathbb{R}(b,g))]\rightarrow W(g)\Big]\Big) \ {\rm by \ transitive}\\
&\geq&\bigwedge_{g\in X}\Big(\mathbb{U}(g)\circledast \Big[\mathbb{R}(d,g)\rightarrow W(g)\Big]\Big)\\
&=&\bigwedge_{g\in X}\Big(\mathbb{U}(g)\circledast \Big[\neg \mathbb{L}(\mathbb{U}_{X-\{d\}})(g)\rightarrow W(g)\Big]\Big)=\mathbb{L}^{\sim 1}(W)(d).\end{eqnarray*}

$\Longleftarrow$. By Remark \ref{remgerl}, $\mathbb{L} =\underline{\mathbb{R}}$ for an {\it Lvr} $\mathbb{R}$. Then  $\forall b, g\in X$,  \begin{eqnarray*}\neg \mathbb{R}(g,b)&=& \underline{\mathbb{R}}^{\sim 1}(\mathbb{U}_{X-\{b\}})(g)\\
&\stackrel{\rm (L4)}{\leq}& \underline{\mathbb{R}}^{\sim 1}\underline{\mathbb{R}} ^{\sim 1}(\mathbb{U}_{X-\{b\}})(g)\\
&=&\bigwedge_{z\in X}\Big(\mathbb{U}(z)\circledast \Big[\neg\underline{\mathbb{R}}(\mathbb{U}_{X-\{g\}})(z)\rightarrow \underline{\mathbb{R}} ^{\sim 1} (\mathbb{U}_{X-\{b\}})(z)\Big]\Big)   \ {\rm by \ Proposition\   \ref{lemmaOIL}(5)}\\
&=&\bigwedge_{z\in X}\Big(\mathbb{U}(z)\circledast \Big[\mathbb{R}(g,z)\rightarrow \neg\mathbb{R}(z,b)\Big]\Big)\\
&=&\bigwedge_{z\in X}\Big(\mathbb{U}(z)\wedge \Big[(\mathbb{U}(z)\rightarrow \mathbb{R}(g,z))\rightarrow \neg\mathbb{R}(z,b)\Big]\Big),\end{eqnarray*} which means that $\forall z\in X$,
\begin{eqnarray*}& &\neg \mathbb{R}(g,b)\leq (\mathbb{U}(z)\rightarrow \mathbb{R}(g,z))\rightarrow \neg\mathbb{R}(z,b)\\
&\Longrightarrow& (\mathbb{U}(z)\rightarrow \mathbb{R}(g,z))\leq \neg \mathbb{R}(g,b)\rightarrow \neg\mathbb{R}(z,b)\\
&\Longrightarrow& (\mathbb{U}(z)\rightarrow \mathbb{R}(g,z))\leq  \mathbb{R}(z,b)\rightarrow \mathbb{R}(g,b)\\
&\Longrightarrow& \mathbb{R}(z,b) \circledast (\mathbb{U}(z)\rightarrow \mathbb{R}(g,z))\leq \mathbb{R}(g,b)\\
&\Longrightarrow& \mathbb{R}(g,z)\circledast (\mathbb{U}(z)\rightarrow \mathbb{R}(z,b))\leq \mathbb{R}(g,b),\end{eqnarray*}
hence $\mathbb{R}$ is transitive.
\end{proof}

\begin{theorem} \label{osutheoremtl}  $\mathbb{L}$ is  a transitive  {\it LVLRAO} iff it satisfies (LT).

(LT) $\forall M, Q\in P(\mathbb{U})$, $ \mathbb{O}\Big(Q,\mathbb{L}(M)\Big)=\mathbb{O}\Big(M, \mathbb{L}^{\sim 1}\mathbb{L}^{\sim 1}(Q)\wedge \mathbb{L}^{\sim 1}(Q)\Big).$
\end{theorem}

\begin{proof} From Theorem \ref{sutheorembgl} and Proposition \ref{oproptl} we only need examine (L)+(L4) $\Longleftrightarrow$ (LT).

(1) (L)+(L4) $\Longrightarrow$ (LT). Indeed $$\mathbb{O}\Big(M, \mathbb{L}^{\sim 1}(Q)\wedge \mathbb{L}^{\sim 1}\mathbb{L}^{\sim 1}(Q)\Big)\stackrel{\rm (L4)}{=}\mathbb{O}\Big(M, \mathbb{L}^{\sim 1}(Q)\Big)\stackrel{\rm (L)}{=}\mathbb{O}\Big(Q,\mathbb{L}(M)\Big),$$ hence (LT) holds.

(1) (LT) $\Longrightarrow$ (L)+(L4). $\forall u\in X$,  take $M=\mathbb{U}_{X-\{u\}}$ in (LT) one get $$\mathbb{L}^{\sim 1}(Q)(u)\stackrel{\rm Proposition\   \ref{lemmaOIL}(2)}{=}\mathbb{O}\Big(Q,\mathbb{L}(\mathbb{U}_{X-\{u\}})\Big)\stackrel{\rm (LT)}{=}\mathbb{O}\Big(\mathbb{U}_{X-\{u\}}, \mathbb{L}^{\sim 1}\mathbb{L}^{\sim 1}(Q)\wedge \mathbb{L}^{\sim 1}(Q)\Big)\leq  \mathbb{L}^{\sim 1}\mathbb{L}^{\sim 1}(Q)(u),$$ so (L4) gains. Using (L4) to (LT) one get (L).
\end{proof}

\begin{proposition} \label{opropsl}    $\mathbb{L}$ is a symmetric LVLRAO iff it fulfills (L1), (L2) and (L5).

(L5) $\forall W\in P(\mathbb{U})$, $\mathbb{L}^{\sim 1}(W)=\mathbb{L}(W)$.
\end{proposition}

\begin{proof}
$\Longrightarrow$. Let $\mathbb{L} =\underline{\mathbb{R}}$ for a symmetric {\it Lvr} $\mathbb{R}$ on $\mathbb{U}$. We only examine (L5). Indeed, $\forall W\in  P(\mathbb{U})$, $b\in X$,
 \begin{eqnarray*}\mathbb{L}^{\sim 1}(W)(b)&=&\bigwedge_{g\in X}\Big(\mathbb{U}(g)\circledast \Big[\neg\underline{\mathbb{R}}(\mathbb{U}_{X-\{b\}})(g)\rightarrow W(g)\Big]\Big)\ {\rm by}\ \mathbb{U}(g)=\mathbb{U}(b)\\
&=&\bigwedge_{g\in X}\Big(\mathbb{U}(b)\circledast \Big[\mathbb{R}(b,g)\rightarrow W(g)\Big]\Big) \ {\rm by \ symmetry}\\
&=&\bigwedge_{g\in X}\Big(\mathbb{U}(b)\circledast \Big[\mathbb{R}(g,b)\rightarrow W(g)\Big]\Big)=\mathbb{L}(W)(b).\end{eqnarray*}

$\Longleftarrow$. By Remark \ref{remgerl}, $\mathbb{L} =\underline{\mathbb{R}}$ for an {\it Lvr} $\mathbb{R}$. Then $\forall b,h\in X$,  \begin{eqnarray*}\mathbb{R}(h,b)=\neg\underline{\mathbb{R}}(\mathbb{U}_{X-\{h\}})(b)\stackrel{\rm (L5)}{=}\neg\underline{\mathbb{R}}^{\sim 1}(\mathbb{U}_{X-\{h\}})(b)\stackrel{\rm Proposition \ref{lemmaOIL}(5)}{=}\neg\underline{\mathbb{R}}(\mathbb{U}_{X-\{b\}})(h)=\mathbb{R}(b,h).\end{eqnarray*} So, $\mathbb{R}$ is symmetric.
\end{proof}

\begin{theorem} \label{osutheoremsl}  $\mathbb{L}$ is  a symmetric   {\it LVLRAO} iff it satisfies (LS).

(LS) $\forall M, Q\in P(\mathbb{U})$, $\mathbb{O}\Big(Q,\mathbb{L}(M)\Big)=\mathbb{O}\Big(M,  \mathbb{L}(Q)\Big).$
\end{theorem}

\begin{proof} From Theorem \ref{sutheorembgl} and Proposition \ref{opropsl} we only need examine (L)+(L5) $\Longleftrightarrow$ (LS).

(1) (L)+(L5) $\Longrightarrow$ (LS).  Indeed $$\mathbb{O}\Big(M, \mathbb{L}(Q)\Big)\stackrel{\rm (L5)}{=}\mathbb{O}\Big(M, \mathbb{L}^{\sim 1}(Q)\Big)\stackrel{\rm (L)}{=}\mathbb{O}\Big(Q,\mathbb{L}(M)\Big).$$

(2) (LS) $\Longrightarrow$ (L)+(L5). $\forall u\in X$,  put $M=\mathbb{U}_{X-\{u\}}$ in (LS) we get $$\mathbb{L}^{\sim 1}(Q)(u)\stackrel{\rm Proposition\   \ref{lemmaOIL}(2)}{=}\mathbb{O}\Big(Q,\mathbb{L}(\mathbb{U}_{X-\{u\}})\Big)\stackrel{\rm (LS)}{=}\mathbb{O}\Big(\mathbb{U}_{X-\{u\}}, \mathbb{L}(Q)\Big)= \mathbb{L}(Q)(u),$$ so (L5) gains. Using (L5) to (LS) one get (L).
\end{proof}

\begin{proposition} \label{opropeul}  $\mathbb{L}$ is an  Euclidean LVLRAO iff it fulfills (L1), (L2) and (L6).

(L6) $\forall W \in P(\mathbb{U})$, $\mathbb{L}\mathbb{L}^{\sim 1} (W)\geq \mathbb{L}^{\sim 1} (W)$.
\end{proposition}

\begin{proof} $\Longrightarrow$. Let $\mathbb{L} =\underline{\mathbb{R}}$ for an  Euclidean {\it Lvr} $\mathbb{R}$ on $\mathbb{U}$. We only examine (L6).  $\forall W\in  P(\mathbb{U})$, $c\in X$,

 \begin{eqnarray*}\mathbb{L}\mathbb{L} ^{\sim 1}(W)(c)&=&\bigwedge_{b\in X}\Big(\mathbb{U}(c)\circledast \Big[\mathbb{R}(b,c)\rightarrow \underline{\mathbb{R}}^{\sim 1}(W)(b)\Big]\Big)\\
&=&\bigwedge_{b\in X}\Big(\mathbb{U}(c)\circledast \Big[\mathbb{R}(b,c)\rightarrow \bigwedge_{g\in X}\Big\{\mathbb{U}(g)\circledast[\neg\underline{\mathbb{R}}(\mathbb{U}_{X-\{b\}})(g)\rightarrow W(g)]\Big\}\Big]\Big)\ {\rm by}\ \mathbb{U}(c)=\mathbb{U}(g)\\
&=&\bigwedge_{b,g\in X}\Big(\mathbb{U}(g)\circledast \Big[\mathbb{R}(b,c)\rightarrow \Big\{\mathbb{U}(g)\wedge [(\mathbb{U}(g)\rightarrow \mathbb{R}(b,g))\rightarrow W(g)]\Big\}\Big]\Big)\ {\rm by}\ \mathbb{U}(g)=\mathbb{U}(b)\\
&=&\bigwedge_{b,g\in X}\Big(\mathbb{U}(g)\circledast \Big[\mathbb{R}(b,c)\rightarrow \Big\{ (\mathbb{U}(b)\rightarrow \mathbb{R}(b,g))\rightarrow W(g)\Big\}\Big]\Big)\\
&=&\bigwedge_{b,g\in X}\Big(\mathbb{U}(g)\circledast \Big[\Big\{ \mathbb{R}(b,c)\circledast (\mathbb{U}(b)\rightarrow \mathbb{R}(b,g))\Big\}\rightarrow W(g)\Big]\Big) \ {\rm by \ Euclidean}\\
&\geq&\bigwedge_{g\in X}\Big(\mathbb{U}(g)\circledast \Big[\mathbb{R}(c,g)\rightarrow W(g)\Big]\Big)\\
&=&\bigwedge_{g\in X}\Big(\mathbb{U}(g)\circledast \Big[\neg\underline{\mathbb{R}}(\mathbb{U}_{X-\{c\}})(g)\rightarrow W(g)\Big]\Big)\\
&=&\mathbb{L}^{\sim 1}(W)(c).\end{eqnarray*}

$\Longleftarrow$. By Remark \ref{remgerl}, $\mathbb{L} =\underline{\mathbb{R}}$ for an {\it Lvr} $\mathbb{R}$. Then $\forall g, z\in X$, \begin{eqnarray*}\neg\mathbb{R}(z,g)&=& \underline{\mathbb{R}}^{\sim 1}(\mathbb{U}_{X-\{g\}})(z)\stackrel{\rm (L6)}{\leq} \underline{\mathbb{R}}\underline{\mathbb{R}}^{\sim 1}(\mathbb{U}_{X-\{g\}})(z)\\
&=&\bigwedge_{b\in X}\Big(\mathbb{U}(z)\circledast \Big[\mathbb{R}(b,z)\rightarrow \underline{\mathbb{R}}^{\sim 1} (\mathbb{U}_{\{g\}})(b)\Big]\Big)\\
&=&\bigwedge_{b\in X}\Big(\mathbb{U}(z)\circledast \Big[\mathbb{R}(b,z)\rightarrow \neg\mathbb{R}(b,g)\Big]\Big)\\
&=&\bigwedge_{b\in X}\Big(\mathbb{U}(z)\wedge \Big[(\mathbb{U}(z)\rightarrow \mathbb{R}(b,z))\rightarrow \neg\mathbb{R}(b,g)\Big]\Big),\end{eqnarray*} which means that $\forall b\in X$,
\begin{eqnarray*}& &\neg \mathbb{R}(z,g)\leq (\mathbb{U}(z)\rightarrow \mathbb{R}(b,z))\rightarrow \neg\mathbb{R}(b,g)\\
&\Longrightarrow& (\mathbb{U}(z)\rightarrow \mathbb{R}(b,z))\leq \neg \mathbb{R}(z,g)\rightarrow \neg\mathbb{R}(b,g)\\
&\Longrightarrow& (\mathbb{U}(b)\rightarrow \mathbb{R}(b,z))\leq  \mathbb{R}(b,g)\rightarrow \mathbb{R}(z,g)\\
&\Longrightarrow& \mathbb{R}(b,g) \circledast (\mathbb{U}(b)\rightarrow \mathbb{R}(b,z))\leq \mathbb{R}(z,g)\\
&\Longrightarrow& \mathbb{R}(b,z)\circledast (\mathbb{U}(b)\rightarrow \mathbb{R}(b,g))\leq \mathbb{R}(z,g),\end{eqnarray*}
i.e., $\mathbb{R}$ is Euclidean.
\end{proof}

\begin{theorem} \label{sutheoremeul}  $\mathbb{L}$ is  an Euclidean  {\it LVLRAO} iff it fulfills  (LE).

(LE) $\forall M, Q\in P(\mathbb{U})$, $\mathbb{O}\Big(Q,\mathbb{L}(M)\Big)=\mathbb{O}\Big(M, \mathbb{L}\mathbb{L}^{\sim 1}(Q)\wedge \mathbb{L}^{\sim 1}(Q)\Big).$
\end{theorem}

\begin{proof} From Theorem \ref{sutheorembgl} and Proposition \ref{opropeul} we only need examine (L)+(L6) $\Longleftrightarrow$ (LE).

(1) (L)+(L6) $\Longrightarrow$ (LE).  Indeed $$\mathbb{O}\Big(M, \mathbb{L}^{\sim 1}(Q)\wedge \mathbb{L}\mathbb{L}^{\sim 1}(Q)\Big)\stackrel{\rm (L6)}{=}\mathbb{O}\Big(M, \mathbb{L}^{\sim 1}(Q)\Big)\stackrel{\rm (L)}{=}\mathbb{O}\Big(Q,\mathbb{L}(M)\Big).$$

(2) (LE) $\Longrightarrow$ (L)+(L6). $\forall o\in X$,  put $M=\mathbb{U}_{X-\{o\}}$ in (LE) we get $$\mathbb{L}^{\sim 1}(Q)(o)=\mathbb{O}\Big(Q,\mathbb{L}(\mathbb{U}_{X-\{o\}})\Big)\stackrel{\rm (LE)}{=}\mathbb{O}\Big(\mathbb{U}_{X-\{o\}}, \mathbb{L}\mathbb{L}^{\sim 1}(Q)\wedge \mathbb{L}^{\sim 1}(Q)\Big)\leq  \mathbb{L}\mathbb{L}^{\sim 1}(Q)(o),$$so  (L6) gains. Using (L6) to (LE) one get (L).
\end{proof}

\begin{proposition} \label{opropml}  $\mathbb{L}$ is a mediate LVLRAO iff it fulfills (L1), (L2) and (L7).

 (L7) $\forall W \in P(\mathbb{U})$, $\mathbb{L}^{\sim 1}\mathbb{L}^{\sim 1}(W)\leq \mathbb{L}^{\sim 1} (W)$.
\end{proposition}
\begin{proof} The proof is similar to Proposition \ref{oproptl}.
\end{proof}

\begin{theorem} \label{osutheoremml}  $\mathbb{L}$ is  a mediate  {\it LVLRAO} iff it satisfies (LM).

(LM) $\forall M, Q\in P(\mathbb{U})$, $ \mathbb{O}\Big(Q,\mathbb{L}(M)\Big)=\mathbb{O}\Big(M, \mathbb{L}^{\sim 1}\mathbb{L}^{\sim 1}(Q)\vee \mathbb{L}^{\sim 1}(Q)\Big).$
\end{theorem}

\begin{proof} The detailed proof is similar to Theorem \ref{osutheoremtl}.
\end{proof}

%
%
%

\subsection{The single axiom characterizations about composite LVLRAOs}

This subsection  provides the single axiom characterizations about LVHRAOs derived from  composition about reflexive, symmetric, transitive, Euclidean $\&$ mediate {\it Lvr}s.

We first consider the combinations of two conditions.

\begin{theorem} \label{osutheoremrtl} (1) $\mathbb{L}$ is a reflexive and transitive, i.e., preordered  {\it LVLRAO} iff it fulfills (LRT).

(LRT) $\forall Q, N\in P(\mathbb{U})$, $\mathbb{O}\Big(Q,\mathbb{L}(N)\Big)=\mathbb{O}\Big(N, Q\wedge \mathbb{L}^{\sim 1}\mathbb{L}^{\sim 1}(Q)\wedge \mathbb{L}^{\sim 1}(Q)\Big).$

(2) $\mathbb{L}$ is a reflexive and symmetric, i.e., tolerance  {\it LVLRAO} iff it fulfills (LRS).

(LRS) $\forall Q, N\in P(\mathbb{U})$, $\mathbb{O}\Big(Q,\mathbb{L}(N)\Big)=\mathbb{O}\Big(N, Q\wedge  \mathbb{L}(Q)\Big).$
\end{theorem}

\begin{proof}  (1) By Theorem \ref{osutheoremrl} and Proposition \ref{oproptl}, we only need examine (LR)+(L4) $\Longleftrightarrow$ (LRT).

(i) (LR)+(L4) $\Longrightarrow$ (LRT). Indeed  $$\mathbb{O}\Big(N, Q\wedge \mathbb{L}^{\sim 1}\mathbb{L}^{\sim 1}(Q)\wedge \mathbb{L}^{\sim 1}(Q)\Big)\stackrel{\rm (L4)}{=}\mathbb{O}\Big(N, Q\wedge \mathbb{L}^{\sim 1}(Q)\Big)\stackrel{\rm (LR)}{=}\mathbb{O}\Big(Q,\mathbb{L}(N)\Big),$$ i.e., (LRT) holds.

(ii) (LRT) $\Longrightarrow$ (LR)+(L4). $\forall b\in X$,  put $N=\mathbb{U}_{X-\{b\}}$ in (LRT) then
$$\mathbb{L}^{\sim 1}(Q)(b)=\mathbb{O}\Big(Q,\mathbb{L}(\mathbb{U}_{X-\{b\}})\Big)\stackrel{\rm (LRT)}{=}\mathbb{O}\Big(\mathbb{U}_{X-\{b\}}, Q\wedge  \mathbb{L}^{\sim 1}\mathbb{L}^{\sim 1}(Q)\wedge \mathbb{L}^{\sim 1}(Q)\Big)\leq  \mathbb{L}^{\sim 1}\mathbb{L}^{\sim 1}(Q)(b),$$ so  (L4) gains. Using (L4) to (LRT) one get (LR).

(2) By Theorem \ref{osutheoremrl} and Proposition \ref{opropsl}, we only need examine (LR)+(L5) $\Longleftrightarrow$ (LRS). The detailed proof is similar to (1).
\end{proof}

\begin{proposition} \label{osutheoremtsl}  (1) $\mathbb{L}$ is a reflexive and  Euclidean  {\it LVLRAO} iff it fulfills (LRE).

(LRE) $\forall Q, N\in P(\mathbb{U})$, $\mathbb{O}\Big(Q,\mathbb{L}(N)\Big)=\mathbb{O}\Big(N, Q\wedge \mathbb{L}\mathbb{L}^{\sim 1}(Q)\wedge \mathbb{L}^{\sim 1}(Q)\Big).$

(2) $\mathbb{L}$ is  a transitive and symmetric  {\it LVLRAO} iff it satisfies (LTS).

(LTS) $\forall Q, N\in P(\mathbb{U})$, $ \mathbb{O}\Big(Q,\mathbb{L}(N)\Big)=\mathbb{O}\Big(N, \mathbb{L}\mathbb{L}(Q)\wedge \mathbb{L}(Q)\Big).$

(3) $\mathbb{L}$ is  a transitive and Euclidean  {\it LVLRAO} iff it satisfies (LTE).

(LTE) $\forall Q, N\in P(\mathbb{U})$, $ \mathbb{O}\Big(Q,\mathbb{L}(N)\Big)=\mathbb{O}\Big(N, \mathbb{L}^{\sim 1}\mathbb{L}^{\sim 1}(Q)\wedge \mathbb{L}\mathbb{L}^{\sim 1}(Q)\wedge \mathbb{L}^{\sim 1}(Q)\Big).$

(4) $\mathbb{L}$ is  a transitive and mediate  {\it LVLRAO} iff it satisfies (LTM).

(LTM) $\forall Q, N\in P(\mathbb{U})$, $ \mathbb{O}\Big(Q,\mathbb{L}(N)\Big)=\mathbb{O}\Big(N, \mathbb{L}^{\sim 1}\mathbb{L}^{\sim 1}(Q)\Big).$

(5)  $\mathbb{L}$ is  a symmetric and Euclidean  {\it LVLRAO} iff it satisfies (LTS).

(6) $\mathbb{L}$ is  a symmetric and mediate  {\it LVLRAO} iff it satisfies (LSM).

(LSM) $\forall Q, N\in P(\mathbb{U})$, $ \mathbb{O}\Big(Q,\mathbb{L}(N)\Big)=\mathbb{O}\Big(N, \mathbb{L}\mathbb{L}(Q)\vee \mathbb{L}(Q)\Big).$
\end{proposition}

\begin{proof} (1) By Theorem \ref{osutheoremrl} and Proposition \ref{opropeul}, we only need examine (LR)+(L6) $\Longleftrightarrow$ (LRE). The detailed proof is similar to Theorem \ref{osutheoremrtl} (1).

(2) From Theorem \ref{osutheoremtl} and Proposition \ref{opropsl} we only need examine (LT)+(L5) $\Longleftrightarrow$ (LTS).

(i) (LT)+(L5) $\Longrightarrow$ (LTS). Indeed $$\mathbb{O}\Big(N, \mathbb{L}(Q)\wedge \mathbb{L}\mathbb{L}(Q)\Big)\stackrel{\rm (L5)}{=}\mathbb{O}\Big(N, \mathbb{L}^{\sim 1}(Q)\wedge \mathbb{L}^{\sim 1}\mathbb{L}^{\sim 1}(Q)\Big)\stackrel{\rm (LT)}{=}\mathbb{O}\Big(Q,\mathbb{L}(N)\Big),$$ i.e., (LTS) holds.

(ii) (LTS) $\Longrightarrow$ (LT)+(L5). $\forall b\in X$,  put $N=\mathbb{U}_{X-\{b\}}$ in (LTS) one have $$\mathbb{L}^{\sim 1}(Q)(b)=\mathbb{O}\Big(Q,\mathbb{L}(\mathbb{U}_{X-\{b\}})\Big)\stackrel{\rm (LTS)}{=}\mathbb{O}\Big(\mathbb{U}_{X-\{b\}}, \mathbb{L}\mathbb{L}(Q)\wedge \mathbb{L}(Q)\Big)\leq  \mathbb{L}(Q)(b).$$  Similarly, one can verify that  $ \mathbb{L}^{\sim 1}(Q)(b)\geq  \mathbb{L}(Q)(b)$. Hence   (L5) gains. Using (L5) to (LTS) one get (LT).

(3) From Theorem \ref{osutheoremtl} and Proposition \ref{opropeul} we only need examine (LT)+(L6) $\Longleftrightarrow$ (LTE). The detailed proof is similar to   (2).

(4) From Theorem \ref{osutheoremtl} and Proposition \ref{opropml} we only need examine (LT)+(L7) $\Longleftrightarrow$ (LTM). The detailed proof is similar to   (2).

(5) It follows by Proposition \ref{osutheoremtsl} (2) and that if $\mathbb{R}$ is symmetric, then $\mathbb{R}$ is transitive iff $\mathbb{R}$ is Euclidean.

(6) From Theorem \ref{osutheoremml} and Proposition \ref{opropsl} we only need examine (LM)+(L5) $\Longleftrightarrow$ (LSM).

(i) (LM)+(L5) $\Longrightarrow$ (LSM). Indeed $$\mathbb{O}\Big(N, \mathbb{L}(Q)\vee \mathbb{L}\mathbb{L}(Q)\Big)\stackrel{\rm (L5)}{=}\mathbb{O}\Big(N, \mathbb{L}^{\sim 1}(Q)\vee \mathbb{L}^{\sim 1}\mathbb{L}^{\sim 1}(Q)\Big)\stackrel{\rm (LM)}{=}\mathbb{O}\Big(Q,\mathbb{L}(N)\Big),$$ i.e., (LSM) holds.

(ii) (LSM) $\Longrightarrow$ (LM)+(L5). $\forall b\in X$,  put $N=\mathbb{U}_{X-\{b\}}$ in (LSM) one have $$\mathbb{L}^{\sim 1}(Q)(b)=\mathbb{O}\Big(Q,\mathbb{L}(\mathbb{U}_{X-\{b\}})\Big)\stackrel{\rm (LSM)}{=}\mathbb{O}\Big(\mathbb{U}_{X-\{b\}}, \mathbb{L}\mathbb{L}(Q)\vee \mathbb{L}(Q)\Big)\geq  \mathbb{L}(Q)(b),$$ similarly, one can verify  $ \mathbb{L}^{\sim 1}(Q)(b)\leq  \mathbb{L}(Q)(b)$. Hence  (L5) gains. Using (L5) to (LSM) one get (LM).
\end{proof}

Next, we consider the combinations of three conditions.

\begin{theorem} \label{osutheoremrstl} $\mathbb{L}$ is a reflexive, transitive and symmetric, i.e., equivalent  {\it LVLRAO} iff it fulfills (LRTS).

(LRTS) $\forall Q, N\in P(\mathbb{U})$,
$\mathbb{O}\Big(Q,\mathbb{L}(N)\Big)=\mathbb{O}\Big(N,  Q\wedge \mathbb{L}(Q)\wedge \mathbb{L}\mathbb{L}(Q)\Big).$
\end{theorem}

\begin{proof} From Theorem \ref{osutheoremrtl} and Proposition \ref{opropsl} we only need examine (LRT)+(L5) $\Longleftrightarrow$ (LRTS).

(1) (LRT)+(L5) $\Longrightarrow$ (LRTS). Indeed $$\mathbb{O}\Big(N,  Q\wedge \mathbb{L}(Q)\wedge \mathbb{L}\mathbb{L}(Q)\Big)\stackrel{\rm (L5)}{=}\mathbb{O}\Big(N, Q\wedge \mathbb{L}^{\sim 1}(Q)\wedge \mathbb{L}^{\sim 1}\mathbb{L}^{\sim 1}(Q)\Big)\stackrel{\rm (LRT)}{=}\mathbb{O}\Big(Q,\mathbb{L}(N)\Big),$$ i.e., (LRTS) holds.

(2) (LRTS) $\Longrightarrow$ (LRT)+(L5). $\forall b\in X$,  put $N=\mathbb{U}_{X-\{b\}}$ in (LRTS) one have $$\mathbb{L}^{\sim 1}(Q)(b)=\mathbb{O}\Big(Q,\mathbb{L}(\mathbb{U}_{X-\{b\}})\Big)\stackrel{\rm (LRTS)}{=}\mathbb{O}\Big(\mathbb{U}_{X-\{b\}}, Q\wedge \mathbb{L}\mathbb{L}(Q)\wedge \mathbb{L}(Q)\Big)\leq  \mathbb{L}(Q)(b),$$ similarly, one obtain   $ \mathbb{L}^{\sim 1}(Q)(b)\geq  \mathbb{L}(Q)(b)$. Hence   (L5) gains. Using (L5) to (LRTS) one get (LRT).
\end{proof}

\begin{proposition} \label{osutheoremrtel} (1) $\mathbb{L}$ is a reflexive, transitive and Euclidean {\it LVLRAO} iff it fulfills (LRTE).

(LRTE) $\forall Q, N\in P(\mathbb{U})$, $\mathbb{O}\Big(Q,\mathbb{L}(N)\Big)=\mathbb{O}\Big(N, Q\wedge \mathbb{L}^{\sim 1}\mathbb{L}^{\sim 1}(Q)\wedge \mathbb{L}\mathbb{L}^{\sim 1}(Q)\wedge \mathbb{L}^{\sim 1}(Q)\Big).$

(2) $\mathbb{L}$ is a reflexive, symmetric and Euclidean  {\it LVLRAO} iff it fulfills (LRTS).

(3) $\mathbb{L}$ is a transitive, symmetric and Euclidean  {\it LVLRAO} iff it fulfills (LTS).
\end{proposition}
\begin{proof} (1) By Theorem \ref{osutheoremrtl} and Proposition \ref{opropeul}, we only need examine (LRT)+(L6) $\Longleftrightarrow$ (LRTE). The detailed proof is similar to Theorem \ref{osutheoremrstl}.

(2) It follows by Theorem \ref{osutheoremrstl} and that if $\mathbb{R}$ is symmetric, then $\mathbb{R}$ is transitive iff $\mathbb{R}$ is Euclidean.

(3) It follows by Proposition \ref{osutheoremtsl} (2) and that if $\mathbb{R}$ is symmetric, then $\mathbb{R}$ is transitive iff $\mathbb{R}$ is Euclidean.
\end{proof}

At last, we give the least (resp., largest) equivalent {\it LVLRAO} on $\mathbb{U}$ and the corresponding {\it Lvr}.

\begin{example} (1) Let $\mathbb{L}_0: P(\mathbb{U})\longrightarrow P(\mathbb{U})$ be defined by $\forall Q\in P(\mathbb{U}), \mathbb{L}_0(Q)=Q$. Then $\forall Q, N\in P(\mathbb{U})$,
$$\mathbb{I}\Big(Q,\mathbb{L}_0(N)\Big)=\mathbb{I}\Big(Q,N\Big)=\mathbb{I}\Big(N,Q\Big)=\mathbb{I}\Big(N,  Q\wedge \mathbb{L}_0(Q)\wedge \mathbb{L}_0\mathbb{L}_0(Q)\Big),$$ i.e., (LRTS) holds.  According to Theorem \ref{osutheoremrstl}, there is an uniquely  equivalent {\it Lvr} $\mathbb{R}_0$ on $\mathbb{U}$  s.t. $\underline{\mathbb{R}_0}=\mathbb{L}_0$. Hence, $\forall a,b \in X$, $$\mathbb{R}_0(a,b)=\neg \underline{\mathbb{R}_0}(\mathbb{U}_{X-\{a\}})(b)=\neg \mathbb{L}_0(\mathbb{U}_{X-\{a\}})(b)=\neg \mathbb{U}_{X-\{a\}}(b)=\mathbb{U}_{\{a\}}(b),$$ which means $\mathbb{R}_0$ is the least equivalent {\it Lvr} on $\mathbb{U}$, and $\mathbb{L}_0$ is the largest equivalent {\it LVLRAO} on $\mathbb{U}$.

(2)  Let $\mathbb{L}_1: P(\mathbb{U})\longrightarrow P(\mathbb{U})$ be defined by $$\forall Q\in P(\mathbb{U}), \forall a\in X,  \mathbb{L}_1(Q)(a)=\mathbb{U}(a)\wedge\bigwedge_{d\in X}Q(d).$$ Then $\mathbb{L}_1(Q)(a)\leq \mathbb{U}(a)\wedge Q(a)=Q(a),$ and
\begin{eqnarray*}\mathbb{L}_1\mathbb{L}_1(Q)(a)&=&\mathbb{U}(a)\wedge\bigwedge_{d\in X}\Big( \mathbb{U}(d) \wedge\bigwedge_{e\in X} Q(e)\Big), \mathbb{U}(a)=\mathbb{U}(d)\\
&=&\mathbb{U}(a)\wedge\bigwedge_{e\in X} Q(e)=\mathbb{L}_1(Q)(a).\end{eqnarray*}
It follows that $\forall Q, N\in P(\mathbb{U})$, \begin{eqnarray*}\mathbb{I}\Big(N,  Q\wedge \mathbb{L}_1(Q)\wedge \mathbb{L}_1\mathbb{L}_1(Q)\Big)&=&\mathbb{I}\Big(N,  \mathbb{L}_1(Q)\Big)\\
&=&\bigwedge_{a,b\in X}\Big(\mathbb{U}(a)\circledast (\neg N(a)\rightarrow Q(b))\Big)\\
&=&\bigwedge_{a,b\in X}\Big(\mathbb{U}(a)\circledast ([\mathbb{U}(a)\circledast (N(a)\rightarrow 0)]\rightarrow Q(b))\Big), \ {\rm by \ Lemma \ref{jbmv}(1)}\\
&=&\bigwedge_{a,b\in X}\Big(\mathbb{U}(a)\circledast ([\mathbb{U}(a)\circledast (Q(b)\rightarrow 0)]\rightarrow N(a))\Big), \mathbb{U}(a)=\mathbb{U}(b)\\
&=&\bigwedge_{a,b\in X}\Big(\mathbb{U}(b)\circledast (\neg Q(b)\rightarrow N(a))\Big)\\
&=&\mathbb{I}\Big(Q,\mathbb{L}_1(N)\Big),
\end{eqnarray*} i.e., (LRTS) holds.  According to Theorem \ref{osutheoremrstl}, there is an uniquely  equivalent {\it Lvr} $\mathbb{R}_1$ on $\mathbb{U}$  s.t. $\underline{\mathbb{R}_1}=\mathbb{L}_1$. Hence, $\forall a,b \in X$, $$\mathbb{R}_1(a,b)=\neg \underline{\mathbb{R}_1}(\mathbb{U}_{X-\{a\}})(b)=\neg \mathbb{L}_1(\mathbb{U}_{X-\{a\}})(b)=\mathbb{U}(b)=\mathbb{U}(a)\wedge \mathbb{U}(b),$$ which means $\mathbb{R}_1$ is the largest equivalent {\it Lvr} on $\mathbb{U}$,  and $\mathbb{L}_1$ is the least equivalent {\it LVLRAO} on $\mathbb{U}$.
\end{example}

\begin{remark}  As we trace the evolution from $L$-rough sets  to $L$-valued rough sets, we can find the following two aspects.

(1) The axiomatic forms of Approach 1 and Approach 2 have progressively become more intricate. This trend can be illustrated through the following details.

 Proposition \ref{fuzzy1}  $\hookrightarrow$ Proposition \ref {theoremL1}, Proposition \ref{fuzzy2} $\hookrightarrow$ Proposition \ref{WSC}.

(2) In contrast, the axiomatic form of Approach 3 has remained unchanged, consistently maintaining a concise structure. This trend can be demonstrated through the following details.

Proposition \ref{fuzzy3}  $\hookrightarrow$ Theorem \ref{sutheorembg}, \ref{osutheoremrst}, \ref{sutheorembgl} and \ref{osutheoremrstl}.

Therefore, axiomatization through Approach 3 has always been an essential research content of rough set and its fuzzy extension \cite{Bao18, JIN232, GL13, CY20, WZ16}.
\end{remark}

\section{Conclusion}

$L$-valued rough sets provide a comprehensive theoretical framework for fuzzy rough sets. Axiomatic characterization is a significant research direction within the field of $L$-valued rough sets. Reflecting on the history of axiomatic research related to (fuzzy) rough sets, we observe that there are generally three approaches to achieve this. Recent literature \cite{LF19,BP21,XW22} has successfully completed the axiomatic characterizations of $L$-valued rough sets using Approaches 1 and 2, yet Approach 3 has remained untouched. So, we presented the axiomatic characterization by Approach 3. The main results are summarized as follows.

(1) We introduced the concepts of inner product (i.e., the degree of intersection), subsethood degree and outer product for $L$-subsets in an $L$-universe, which are extensions of the corresponding concepts for  $L$-subsets within  a classical universe. These concepts  provide feasible tool for the study of fuzzy rough set, fuzzy order and fuzzy topology within the $L$-universe framework.

(2) By utilizing the defined inner product (outer product), we characterized 18 types of {\it LVHRAO}s ({\it LVLRAO}s) derived from general, reflexive, symmetric, transitive, Euclidean, and mediate {\it Lvr}s, along with their combinations. This makes that the results regarding the single axiom characterizations of classical and fuzzy rough sets in \cite{Bao18,JIN232,GL13,CY20,WZ16} can be viewed as direct consequences of our findings. Additionally, the method presented in this paper also offers a reference point for axiomatic research on other types of $L$-valued rough sets, such as the variable precision $L$-valued rough set based on $L$-valued relations \cite{AT24} and  $L$-valued rough set based on $L$-valued  coverings \cite{EI24}.

This paper presents numerous axiomatic descriptions of $L$-valued rough approximation operators. To facilitate the readers, we   compiled these characterizations into Table \ref{tac} for easy reference.

\begin{table}[!htbp] \footnotesize
\begin{center}\tabcolsep 0.05in\doublerulesep 0in\renewcommand{\arraystretch}{1.2}
\centering\caption{Single axiomatic characterizations about {\it LVHRAO}s and {\it LVLRAO}s} \label{tac}
\begin{tabular}{llllllllll}
\hline
       $L$-valued relation& \ \ \ \ \  \ \ \ \ \ \ \ \ \  \ \ \ \ {\it LVHRAO} & & \ \ \ \ \  \ \ \ \ \ \ \ \ \  \ \ \ \  {\it LVLRAO}& \\

         & Axiom &Position & Axiom  & Position \\
\hline
General & (H) &Theorem \ref{sutheorembg} &  (L) & Theorem \ref{sutheorembgl} \\

Reflexive  & (HR) & Theorem \ref{osutheoremr}& (LR) & Theorem \ref{osutheoremrl}   \\

Transitive  & (HT) & Theorem \ref{osutheoremt} & (LT) & Theorem \ref{osutheoremtl} \\

Symmetric  & (HS) & Theorem \ref{osutheorems} & (LS) & Theorem \ref{osutheoremsl}\\

Euclidean  & (HE) & Theorem \ref{sutheoremeu} & (LE) & Theorem \ref{sutheoremeul}\\

Mediate  & (HM)& Theorem \ref{osutheoremm} & (LM) & Theorem \ref{osutheoremml}\\
\hline
Reflexive and transitive & (HRT)& & (LRT)  \\

 & & Theorem \ref{osutheoremrt} &   & Theorem \ref{osutheoremrtl}\\

Reflexive and symmetric & (HRS) & & (LRS)\\
\hline
Reflexive and Euclidean & (HRE) & & (LRE)  \\

Transitive and symmetric & (HTS)&  & (LTS)\\

Transitive and  Euclidean & (HTE) &Proposition \ref{osutheoremts}& (LTE) & Proposition \ref{osutheoremtsl} \\

Transitive and mediate & (HTM) & &(LTM)  \\

Symmetric and Euclidean & (HTS)&  & (LTS)  \\

Symmetric and mediate & (HSM)  &  &  (LSM) \\
\hline
Reflexive, transitive and symmetric  & (HRTS) & Theorem \ref{osutheoremrst}& (LRTS)  & Theorem \ref{osutheoremrstl} \\
\hline
Reflexive, transitive and Euclidean  & (HRTE)& & (LRTE) \\

Reflexive, symmetric and Euclidean  & (HRTS) &Proposition \ref{osutheoremrte}& (LRTS)   &  Proposition \ref{osutheoremrtel} \\

Transitive, symmetric and Euclidean  & (HTS) &  & (LTS)  \\
 \hline
\end{tabular}
\end{center}
\end{table}

\section*{Acknowledgements}

This work was supported by National Natural Science Foundation of China (No.12171220).

%

\section*{Declaration of competing interest}

The authors declare that they have no known competing financial interests or personal relationships that could have appeared to influence the work reported in this paper.

%
%

\end{document}